\documentclass[final]{siamart1116}
\pdfoutput=1
\newsiamthm{remark}{Remark}

\usepackage{lipsum}
\usepackage{amsfonts}
\usepackage{graphicx}
\usepackage{epstopdf}
\usepackage{algorithmic}
\Crefname{ALC@unique}{Line}{Lines}
\usepackage{graphicx}
\usepackage{epstopdf}
\usepackage{amsmath}
\usepackage{mathtools}
\usepackage{subcaption}

\ifpdf
  \DeclareGraphicsExtensions{.eps,.pdf,.png,.jpg}
\else
  \DeclareGraphicsExtensions{.eps}
\fi

\numberwithin{theorem}{section}

\newcommand{\TheTitle}{An Automated Singularity-Capturing Scheme for \\ Fractional Differential Equations} 
\newcommand{\TheAuthors}{Jorge Suzuki and Mohsen Zayernouri}

\newcommand{\RunTitle}{An Automated Singularity-Capturing Scheme for FDEs} 
\headers{\RunTitle}{\TheAuthors}

\title{{\TheTitle}}

\author{
  Jorge Suzuki \thanks{Department of Computational Mathematics, Science and Engineering, Michigan State University
    (\email{suzukijo@msu.edu}, \email{suzukijo@egr.msu.edu}.}
  \and
  Mohsen Zayernouri \thanks{Department of Computational Mathematics, Science and Engineering, Michigan State University (\email{zayern@msu.edu}, \email{zayern@egr.msu.edu}).}
}

\usepackage{amsopn}

\begin{document}

\maketitle

\begin{abstract}
Solutions to fractional models inherently exhibit non-smooth behavior, which significantly deteriorates the accuracy and therefore efficiency of existing numerical methods. We develop a two-stage data-infused computational framework for accurate time-integration of single- and multi-term fractional differential equations. In the first stage, we formulate a self-singularity-capturing scheme, given available/observable data for diminutive time. In this approach, the fractional differential equation provides the necessary knowledge/insight on how the hidden singularity can bridge between the initial and the subsequent short-time solution data. We develop a new self-singularity-capturing finite-difference algorithm for automatic determination of the underlying power-law singularities nearby the initial data, employing gradient descent optimization. In the second stage, we can utilize the multi-singular behavior of solution in a variety of numerical methods, without resorting to making any \textit{ad-hoc}/uneducated guesses for the solution singularities. Particularly, we employed an implicit finite-difference method, where the captured singularities, in the first stage, are taken into account through some Lubich-like correction terms, leading to an accuracy of order $\mathcal{O}(\Delta t^{3-\alpha})$. Our computational results demonstrate that the developed framework can either fully capture or successfully control the solution error in the time-integration of fractional differential equations, especially in the presence of strong multi-singularities.
\end{abstract}



\section{Introduction}
\label{Sec:Intro}

Fractional differential equations (FDEs) have been successfully applied in diverse problems that present the fingerprint of power-laws/heavy-tailed statistics, such as visco-elastic modeling of bio-tissues \cite{Craiem2008, Magin2010, naghibolhosseini2015estimation, naghibolhosseini2018, Perdikaris2014}, cell rheology behavior \cite{Djordjevic2003}, food preservatives \cite{Jaishankar2013}, complex fluids \cite{Jaishankar2014}, visco-elasto-plastic modeling for power-law-dependent stresses/strains \cite{Suzuki2016, Suzuki2017Thesis, Hei2018}, earth sciences \cite{ZHANG201729}, among others.

The increasing number of works involving fractional modeling in the past two decades is largely due to the development of efficient numerical schemes for fractional-order/partial differential equations (FDEs/FPDEs). Starting with Lubich's early works \cite{Lubich1983, Lubich1986} on finite-difference schemes, and followed by several developments, \textit{e.g.}, on discretizing Burger's equation \cite{Sugimoto1991}, fractional Adams methods for nonlinear problems \cite{Diethelm2002, Diethelm2004, Zayernouri2016MS}, and fractional diffusion processes \cite{Sun2006, Lin2007, Gorenflo2002, Gao2014}. Also, other classes of global methods were developed, such as spectral methods for FDEs/FPDEs \cite{Zayernouri2014FracDelay,Zayernouri_FODEs_2014, Zayernouri15-SIAM-TFSLP, zayernouri2015-Variable,kharazmi2017sem,Samiee2016,samiee2017unified,kharazmi2017FSEM} and for distributed-order differential equations (DODEs) \cite{kharazmi2017petrov,kharazmi2018fractional,samiee2018petrov}, which are efficient for low-to-high dimensional problems. Of particular interest, we highlight a family of fast convolution schemes for time-fractional integrals/derivatives, initially developed by Lubich \cite{Lubich2002}. In such schemes, instead of a direct discretization by splitting the fractional operator in $N$ uniform time-integration intervals, it is split in exponentially increasing intervals, reducing the overall computational complexity from $\mathcal{O}(N^2)$ to $\approx \mathcal{O}(N \log N)$. Furthermore, the power-law kernels are computed through a convolution quadrature using complex integration contours \cite{Lubich2002, Hiptmair2003, Schadle2006, Li2010, Yu2016JCP}. The storage requirements are also reduced from $\mathcal{O}(N)$ to $\approx \mathcal{O}(\log N)$. Recently, Zeng \textit{.et.al} \cite{Zheng2017} developed an improved version of the scheme using a Gauss-Laguerre convolution quadrature, which utilizes real-valued integration contours. In the same work, the authors also addressed the short memory principle and applied Lubich's correction terms \cite{Lubich1986}, leading to a more accurate and stable fast time-stepping scheme. We also refer readers to other classes of fast schemes for time-fractional operators that make use of the resulting matrix structures \cite{Breiten2016} and kernel compression \cite{Baffet2016} techniques.

There have also been works on `tunable accuracy', spectral collocation methods for FDEs/FPDEs that utilize singular basis functions \cite{Zayernouri_FODEs_2014, Zeng2016, Lischke2017}. Zayernouri and Karniadakis \cite{Zayernouri_FODEs_2014} developed an exponentially accurate spectral collocation method utilizing fractional Lagrange basis functions, given by the product of a singular term with real power $\mu$ and a smooth part given by Lagrange interpolands. Zeng \textit{.et.al} \cite{Zeng2016} later generalized the scheme for variable-order FDEs/FPDEs with endpoint singularities and demonstrated that a proper tuning of the power $\mu$ enhanced the accuracy of the numerical solutions. Lischke \textit{.et.al} \cite{Lischke2017} developed a Laguerre Petrov-Galerkin method for multi-term FDEs with tunable accuracy and linear computational complexity. Despite the high precision obtained by fine-tuning the bases in the developed works, such numerical accuracy is extremely sensitive to the 'single'-singularity basis parameter $\mu$ and no self-capturing scheme was developed to find its correct, application-specific value.

The accuracy and efficiency of the aforementioned numerical schemes strongly depend on the regularity of the solution. To address such problem, a correction method was introduced by Lubich \cite{Lubich1986}, which utilizes $M$ correction terms with singular powers $\sigma_k$, with $k=1,\,\dots,\,M$, given regularity assumptions for $u(t)$. The corresponding correction weights are obtained through the solution of a Vandermonde-type linear system of size $M$. Given a convenient choice of singular powers, some developed schemes were able to attain their theoretical accuracy with non-smooth solutions \cite{Zeng2017CMAME, Zheng2017}. However, such works utilize \textit{ad-hoc} choices of $\sigma_k$, even without prior knowledge of the regularity of the solution. According to Zeng \textit{.et.al} \cite{Zheng2017}, for such cases, a reasonable choice would be $\sigma_k = 0.1 k$, which could improve the numerical accuracy when strong singularities are present, and without loss of accuracy when there is enough regularity for $u(t)$. However, such procedure could still lead to arbitrary choices of $\sigma_k$ far from the true singularities of the solution. There are higher chances of approximating the singularities by increasing the number of correction terms, but using $M > 9$ terms induces a large condition number on the Vandermonde system (see Figure \ref{fig:Cond_a}), and consequently large residuals, leading to large errors when computing the correction weights, which are propagated to the discretized fractional operator \cite{Diethelm2006}. 

Regarding the numerical solution of FDEs with correction terms, let $N$ be the total number of time-steps with size $\Delta t$. For the first $M$ time-steps, prior information of the numerical solution for $u(t)$, see (\ref{eq:discr_2}), is required. The procedure proposed in \cite{Zheng2017} involves the solution of a sub-problem with time domain $t \in (0,\,M \Delta t]$ using a time-step size $\tau = \Delta t^2$ and one correction term to obtain the numerical solution $u^N_n$ for $n = 1,\dots, M$. However, the choice of smaller step size $\tau$ even for a short time might be too expensive, and the numerical solution for the $M$ initial time-steps might experience total loss of accuracy in the presence of strong singularities, not captured by the single guessed correction term. Such deficiencies outline the need for \textbf{I)} Capturing the multi-singularities of $u(t)$ and \textbf{II)} Efficient and accurate schemes that incorporate all captured singularities for time-integration of the (most crucial) initial $M$ time-steps.

The main contribution of the present work is to develop a two-stage data-infused computational framework for accurate time-integration of FDEs. In Stage-I, we formulate a self-singularity-capturing framework where: 
\begin{itemize}
	\item The scheme utilizes available data for initial diminutive time (application-oriented), and self-captures/determines multi-singular behaviors of the solution through the knowledge introduced by the FDE and its corresponding fractional operators.
	
	\item We develop a new (finite-difference based) algorithm for automatic determination of the underlying power-law singularities nearby the initial data, employing gradient descent optimization. The singularities are introduced through Lubich's correction method \cite{Lubich1986}.
	
	\item We introduce the capturing scheme for $M$ correction terms and construct a hierarchical, self-capturing framework. We test the framework for the particular case of up to $S = 3$ singularities and $M=3$ correction terms.
	
	\item The self-capturing procedure makes use of two stopping criteria for the error minimization, namely $\epsilon$ (for solution error) and $\epsilon_1$ (for gradient norm error), where the numerical convergence with respect to $\epsilon$ defines if the singularities are captured. Numerical convergence towards the tolerance $\epsilon_1$ defines if additional correction terms are needed to capture the true solution singularities.
\end{itemize}     
In Stage-II we utilize the captured multi-singularities for the full time-integration of the FDE in question, where a variety of the aforementioned numerical methods could be, in general, employed. Our approach is explained as the following:
\begin{itemize}
	\item To handle the numerical solution of the FDE using multiple correction terms, we develop an implicit finite-difference method, where we solve a linear system of size $M$ for the first $M$ time-steps. Therefore, we incorporate the captured singular behavior up to the desired precision $\epsilon$, without the need of using a fine time-grid;
	
	\item We numerically demonstrate that the developed methodology is of order $\mathcal{O}(\Delta t^{3-\alpha})$;	
\end{itemize}    

We perform a set of numerical tests, where we demonstrate that the developed scheme is able to capture multiple singularities using a few number of time-steps, which can later be utilized for efficient time-integration of FDEs with relatively large time-step size $\Delta t$, being much more accurate than \textit{ad-hoc} choices of singularities when the regularity of the solution is unknown. The successful capturing of singularities motivates the development of kernel- and knowledge-based refinement of time-grids near $t = 0$ \cite{Stynes2017}, as well as self-construction of basis function spaces using M\"{u}ntz polynomials \cite{Esmaeili2011, Hou2017} for spectral element methods for FDEs/FPDEs \cite{kharazmi2017FSEM}.

This paper is organized as follows: Section \ref{Sec:Definitions} introduces the definitions for fractional operators and the FDE in consideration. In section \ref{Sec:TwoStage} we develop the two-stage framework for efficient time-integration of FDEs, where we start with the self-singularity-capturing scheme (Stage-I) in Section \ref{Sec:Stage-I}, followed by the finite-difference scheme for solution of FDEs (Stage-II) in Section \ref{Sec:Stage-II} using the captured singularities. The numerical results with discussions for self-capturing up to three singularities are shown in Section \ref{Sec:Numerical}, followed by the conclusions in Section \ref{Sec:Conclusions}.

\section{Definitions}
\label{Sec:Definitions}

We start with some preliminary definitions of fractional calculus \cite{Podlubny99}. The left-sided Riemann-Liouville (RL) integrals of order $\alpha \in \mathbb{R}$, with $0<\alpha<1$ and $t \in \mathbb{R}$, is defined, as
\begin{equation}
\label{Eq: left RL integral}
(\prescript{RL}{a}{\mathcal{I}}_{t}^{\alpha} v) (t) = \frac{1}{\Gamma(\alpha)} \int_{a}^{t} \frac{v(s)}{(t - s)^{1-\alpha} }\, ds,\,\,\,\,\,\, t>a, 
\end{equation}
where $\Gamma$ represents the Euler gamma function and $a$ denotes the lower integration limit. The corresponding inverse operator, i.e., the left-sided RL fractional derivatives of order $\alpha$, is then defined based on (\ref{Eq: left RL integral}) as  
\begin{equation}
\label{Eq: left RL derivative}
(\prescript{RL}{a}{\mathcal{D}}_{t}^{\alpha} v) (t) = \frac{d}{dt} (\prescript{RL}{a}{\mathcal{I}}_{t}^{1-\alpha} v) (t) = \frac{1}{\Gamma(1-\alpha)}  \frac{d}{dt} \int_{a}^{t} \frac{v(s)}{(t - s)^{\alpha} }\, ds,\,\,\,\,\,\, t>t_L.
\end{equation}
Furthermore, the corresponding left-sided Caputo derivatives of order $\mu \in (0,1)$ is obtained as 
\begin{equation}
\label{Eq: left Caputo derivative}
(\prescript{C}{a}{\mathcal{D}}_{t}^{\mu} v) (t) = (\prescript{RL}{a}{\mathcal{I}}_{t}^{1-\alpha} \frac{dv}{dx}) (x) = \frac{1}{\Gamma(1-\alpha)}  \int_{a}^{t} \frac{v^{\prime}(s)}{(t - s)^{\alpha} }\, ds,\,\,\,\,\,\, t>a.
\end{equation}
The definitions of Riemann-Liouville and Caputo derivatives are linked by the following relationship, which can be derived by a direct calculation
\begin{equation}
\label{Eq: Caputo vs. Riemann}
(\prescript{RL}{a}{\mathcal{D}}_{t}^{\alpha} v) (t)  =  \frac{v(a)}{\Gamma(1-\alpha) (t-a)^{\alpha}}  +   (\prescript{C}{a}{\mathcal{D}}_{t}^{\mu} v) (t),
\end{equation}
which denotes that the definition of the aforementioned derivatives coincide when dealing with homogeneous Dirichlet initial/boundary conditions.

We now introduce the Cauchy problem to be solved in this work and its corresponding well-posedness (see Theorem 3.25\textbf{(i)} \cite{Kilbas2006} with $n = 1$). Let $C(\Omega)$ be the space of continuous functions $u(t)$ in $\Omega=[a,b]$ with the norm:
\begin{equation}
\lvert \lvert u \rvert \rvert_{C(\Omega)} = \max_{t \in \Omega}\lvert u(t) \rvert
\end{equation}
Also, let $0 < \alpha < 1$ and $\gamma \in \mathbb{R^+}$, with $\gamma \le \alpha$. We define $C_\gamma(\Omega)$ to be the following weighted space of continuous functions:
\begin{equation}
C_\gamma(\Omega) = \lbrace g(t)\big{\vert}\,(t-a)^\gamma g(t) \in C(\Omega),\, \vert\vert g \vert\vert_{C_\gamma} = \vert\vert(t-a)^\gamma g(t)\vert\vert_C \rbrace,
\end{equation}
Let $G$ be an open set in $\mathbb{R}$, and let a function $f(t,u(t)):(a,b]\times G \to \mathbb{R}$, such that $f(t,u(t)) \in C_\gamma(\Omega)$ and is Lipschitz continuous with respect to any $u(t) \in G$. Let the fractional Cauchy problem of interest given by:
\begin{equation}\label{eq:problemdef1}
\prescript{C}{a}{}\mathcal{D}^\alpha_t u(t) = f(t,u), \quad u(a) = u_0.
\end{equation}
Given the aforementioned conditions, there exists a unique solution for (\ref{eq:problemdef1}) in the following space of functions:
\begin{equation}
\mathbf{C}^\alpha_{\gamma}(\Omega) = \lbrace u(t) \in C(\Omega)\big{\vert}\, \prescript{C}{a}{}\mathcal{D}^\alpha_t u(t) \in C_\gamma(\Omega)\rbrace.
\end{equation}
The corresponding solution with respect to $u(t)$ is equivalent to the following Volterra integral equation of second kind: $u(t) = u_0 + \frac{1}{\Gamma(\alpha)}\int^t_a \frac{f(s,u(s))}{(t-s)^{1-\alpha}}\,ds$, with $a \le t \le b$. We remark that for other particular forms of (\ref{eq:problemdef1}), the corresponding Volterra equation of second kind contains Mittag-Leffler functions $E_\alpha(z)$, which are defined through infinite sums and thus are not computationally-friendly. Therefore, we choose to work with the presented form, since the developed self-singularity-capturing framework can be extended for other FDEs in a straightforward way, without resorting to computationally expensive hyper-geometric functions. 
Furthermore, although the existence and uniqueness of Cauchy problems has been investigated \cite{Samko1993, Delbosco1996, Kilbas2006, Diethelm2002}, yet, there is no comprehensive framework to understand the singular behavior of the solution given any general $f(t,u(t)) \in C_\gamma (\Omega)$. The developed formulation in this work is suitable for problems where a few initial discrete data points are known from $u(t)$ and $f(t,u(t))$.

Due to the employed finite-difference discretization over the RL derivative in this work (see Section \ref{Sec:FDun}), we choose to rewrite (\ref{eq:problemdef1}) using (\ref{Eq: Caputo vs. Riemann}) to obtain:
\begin{equation}\label{eq:problemdef}
\prescript{RL}{a}{}\mathcal{D}^\alpha_t u(t) = f(t,u) - \frac{u_0}{\Gamma(1-\alpha) (t-a)^{\alpha}}.
\end{equation}
As will be shown in the next sections, the developed scheme is able to capture singularities with a minimal number of time-steps, \textit{i.e.} with restricted data for a short period of time.

\section{Two-Stage Time-Integration Framework}
\label{Sec:TwoStage}

In this section, we develop the two-stage framework for efficient time-integration of FDEs. We start with Stage-I, where the singularity-self-capturing scheme is presented for $S$ singularities, computed over a small number of initial data points denoted by $\tilde{N}$. Then, after the multi-singular solution behavior is learned, we present the developed solution for FDEs in Stage-II, where the subsequent full time-integration of the FDE is carried out over $N$ time steps (see Figure \ref{fig:time_grid}).

\subsection{Stage-I: Self-Singularity-Capturing Stage}
\label{Sec:Stage-I}

We develop the self-singu\-larity-capturing framework, starting with $M$ correction terms for the initial time-steps $\tilde{N}$, at which data is available. We then utilize the developed algorithm in a self-capturing approach, through a hierarchical and iterative fashion. 

\subsubsection{Minimization Method}
\label{Sec:Minimization}

\begin{figure}[t!]
	\centering
	\includegraphics[width=0.55\textwidth]{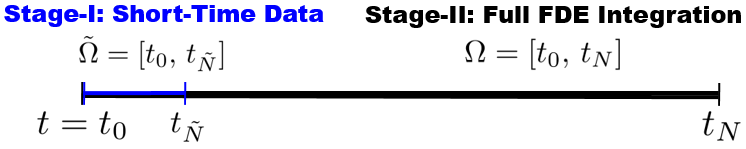}
	\caption{In Stage-I, the scheme captures the singularities for the short time-domain $\tilde{\Omega}$, given the short initial data. In Stage-II, the framework uses the singularities obtained in Stage-I for full time-integration of the FDE over $\Omega$. \label{fig:time_grid}}
\end{figure}

Let $t \in \tilde{\Omega}$, with $\tilde{\Omega} = [t_0,\,t_{\tilde{N}}]$, where $\tilde{N}$ denotes the number of initial short data points, such that $\tilde{N} \ge M$, and let $\boldsymbol{\sigma} \in \mathbb{R}^M$ be the following correction-power (singularity) vector:
\begin{equation}
\boldsymbol{\sigma} = \left[ \sigma_1,\,\, \sigma_2,\,\,\dots,\,\, \sigma_M \right]^T.
\end{equation}
We define an error function $E:\mathbb{R}^M \to \mathbb{R}^+$, given by:
\begin{equation} \label{eq:Cost}
E(\boldsymbol{\sigma}) = \sum^{\tilde{N}}_{n=1} \left( u^{data}_n - u^N_n(\boldsymbol{\sigma}) \right)^2,
\end{equation}
where $u^N_n(\boldsymbol{\sigma})$ denotes the $\boldsymbol{\sigma}$-dependent numerical solution of $u(t)$ at $t = t_n$, and $u^{data}_n$ represents the known initial shot-time data. Figure \ref{fig:time_grid} illustrates the integration domain $\tilde{\Omega}$ for Stage-I, where the error function (\ref{eq:Cost}) is evaluated for a short time.
We apply an iterative gradient descent scheme in order to find $\boldsymbol{\sigma}^*$ that minimizes (\ref{eq:Cost}). Let $\boldsymbol{\sigma}^{k}$ be known at the $k$-th iteration. We compute the updated value $\boldsymbol{\sigma}^{k+1}$ for iteration $k+1$ in the following fashion \cite{Boyd2004}:
\begin{equation}
\boldsymbol{\sigma}^{k+1} = \boldsymbol{\sigma}^k - \gamma^k \nabla E(\boldsymbol{\sigma}^k),
\end{equation}
where $\gamma^k$ denotes the following two-point step size for the $k-$th iteration, given by \cite{Barzilai1988}:
\begin{equation}\label{eq:step}
\gamma^k = \frac{\left(\boldsymbol{\sigma}^k - \boldsymbol{\sigma}^{k-1}\right)^T \left[ \nabla E(\boldsymbol{\sigma}^k) - \nabla E(\boldsymbol{\sigma}^{k-1}) \right]}{\vert\vert \nabla E(\boldsymbol{\sigma}^k) - \nabla E(\boldsymbol{\sigma}^{k-1}) \vert \vert^2_{L^2(\mathbb{R}^M)}}.
\end{equation}
The gradient $\nabla E:\mathbb{R}^+ \to \mathbb{R}^M$ of the error function is given by $\nabla E(\boldsymbol{\sigma}) = [\partial E/\partial \sigma_1, \,\,\allowbreak \partial E/\partial \sigma_2, \,\,\dots,\,\, \partial E/\partial \sigma_M]^T$, where, instead of obtaining a closed form for the derivatives (which would implicitly involve the differentiation of $u^N_i$), we utilize complex-step differentiation. Therefore, let $\mathbf{e}_j = [0,\,\,0,\,\,\dots,\,\,1,\,\,\dots,\,\,0]^T$ be a vector in $\mathbb{R}^M$ of zeros, with a unit value in the $j$-th entry. Therefore, we have the following:
\begin{equation}\label{eq:complex}
\nabla E(\boldsymbol{\sigma})_j \approx \frac{Im(E(\boldsymbol{\sigma} + i \Delta \sigma\mathbf{e}_j))}{\Delta \sigma}, \quad j = 1,\,2,\,\dots,\, M, \quad i=\sqrt{(-1)},
\end{equation}
which only requires the additional evaluation of the error function perturbed by $\Delta \sigma$, which can be taken, \textit{e.g.}, as $10^{-14}$. Given two numerical tolerances $\epsilon$ and $\epsilon_1$, we iterate and find a new $\boldsymbol{\sigma}^{k+1}$ while both criteria $E(\boldsymbol{\sigma}^{k+1}) > \epsilon$ and $||\nabla E(\boldsymbol{\sigma}^k)|| > \epsilon_1$ are true. The latter criterion is introduced to minimize our error function, while the former is used for error control in our self-capturing scheme. As will be shown in the next sections, the criterion $||\nabla E(\boldsymbol{\sigma}^k)|| < \epsilon_1$ is satisfied for small $\epsilon_1$ even when $M < S$, but the criterion $E(\boldsymbol{\sigma}^{k+1}) < \epsilon$ will only be satisfied when we use enough correction terms (see Section \ref{Sec:Self-framework}) to fully capture the $S$ number of power-law singularities.

\subsubsection{Numerical Scheme for Short-Time Integration}
\label{Sec:FDun}


In the singularity-capturing scheme, given $\boldsymbol{\sigma}$, in order to compute the error in (\ref{eq:Cost}) and the gradient (due to the error perturbation) (\ref{eq:complex}), we need to compute the numerical solution $u^N_n$ for $n=1,\, \dots, \,\tilde{N}$ initial data points over time-intervals $\Delta t$. For this purpose, we employ the finite-difference discretization with corrections presented in \cite{Lubich1986} for the fractional RL derivative. Including correction terms, the discretization of the left-sided RL fractional derivative of order $\alpha$, evaluated at time $t = t_n$, with initial time $a = 0$, has the following form:
\begin{equation}
\prescript{RL}{0}{}\mathcal{D}^\alpha_t u(t) \big|_{t = t_n} \approx \prescript{RL}{0}{}\mathcal{D}^{(\alpha, M)}_t u^N(t) \big|_{t = t_n}.
\end{equation}
The above approximation is augmented by the \textit{so-called} Lubich's correction terms, which is given by:
\begin{equation}\label{eq:corrected_deriv}
\prescript{RL}{0}{}\mathcal{D}^{(\alpha, M)}_t u^N(t) \big|_{t = t_n} = \prescript{RL}{0}{}\mathcal{D}^{\alpha}_t u^N(t) \big|_{t = t_n} + \sum^M_{j=1} W_{j,n}(\boldsymbol{\sigma}) \left(u^N_j - u_0 \right),
\end{equation}
with $j=1, \dots, M$. The term $W_{j,n}(\boldsymbol{\sigma})$ denotes the correction weights, which depend on the correction powers $\boldsymbol{\sigma} \in \mathbb{R}^M$. We assume a power-law solution singularity about $t=0$, driven by the power-law kernel of the RL fractional derivative, \textit{i.e.}, $u\big{\vert}_{\substack{nearby \\ t = 0}} \approx \sum^M_{k=1} c_k t^{\sigma_k}$. If $\boldsymbol{\sigma} = \lbrace \sigma_k \rbrace^M_{k=1}$ are known, the RL fractional derivative of each singular term, $t^{\sigma_k}$,  can be obtained by:
\begin{equation}\label{eq:discretized_correction}
	\prescript{RL}{0}{}\mathcal{D}^\alpha_t (t^{\sigma_k}) \big|_{t = t_n} +  \sum^M_{j=1} W_{j,n} t^{\sigma_k}_j = \frac{\Gamma(1+\sigma_k)}{\Gamma(1+\sigma_k - \alpha)} t^{\sigma_k - \alpha}_n, \quad k=1, \dots, M,
\end{equation}
where the first term on the left-hand side of (\ref{eq:discretized_correction}) denotes the discretized RL fractional derivative of $t^{\sigma_k}$, while the right-hand side represents the exact fractional derivative of $t^{\sigma_k}$, evaluated at $t = t_n$. Equation (\ref{eq:discretized_correction}) can be written as:
\begin{equation}\label{eq:final_weights}
	\sum^M_{j=1} V_{k,j} W_{j,n} = \frac{\Gamma(1+\sigma_k)}{\Gamma(1+\sigma_k - \alpha) \Delta t^\alpha} n^{\sigma_k - \alpha} - \frac{1}{\Delta t^{\sigma_k}} \prescript{RL}{0}{}\mathcal{D}^\alpha_t (t^{\sigma_k}) \big|_{t = t_n},
\end{equation}
with $k=1, \dots, M$, where $V_{k,j} = j^{\sigma_k}$ denotes a Vandermonde matrix with size $M \times M$. Therefore, to obtain the starting weights $W_{j,n}$, the above linear system has to be solved for all $n = 1,\, \dots,\, \tilde{N}$. 

Substituting (\ref{eq:corrected_deriv}) into (\ref{eq:problemdef}), and assuming homogeneous initial conditions $u_0 = 0$, we obtain the following discrete form for the FDE:
\begin{equation}\label{eq:discretized_initial}
\prescript{RL}{0}{}\mathcal{D}^{\alpha}_t u^N(t) \big|_{t = t_n} +  \sum^M_{j=1} W_{j,n}(\boldsymbol{\sigma}) \left(u^N_j - u_0 \right) = f^{data}_n, \quad n=1,\, \dots,\,\tilde{N}.
\end{equation}

In order to discretize $\prescript{RL}{0}{}\mathcal{D}^{\alpha}_t u^N(t)$, we follow the difference scheme developed in \cite{Zheng2017}, which is based on a second-order interpolation of $u(t)$ for the fractional RL integral (\ref{Eq: left RL integral}), similar to the L1-2 scheme for fractional Caputo derivatives developed in \cite{Gao2014}. Therefore, when evaluating (\ref{Eq: left RL integral}) at $t=t_n$, we split it into local $\mathcal{L}$ and history $\mathcal{H}$ parts as follows:
\begin{equation}\label{eq:conv_short_long}
	\prescript{RL}{0}{\mathcal{I}}_{t}^{\alpha} u^N(t)\big |_{t=t_n} = \underbrace{\int^{t_n}_{t_{n-1}} k(t_n-s) u^N(s)\,ds}_\text{$\mathcal{L}(u,t_n)$} + \underbrace{ \sum^{n-2}_{k=0} \int^{t_{k+1}}_{t_k} k(t_n-s) u^N(s)\,ds}_\text{$\mathcal{H}(u,t_n)$},
\end{equation}
with the kernel $k(t) = t^{\alpha-1}/\Gamma(\alpha)$. The function $u^N(t)$ is approximated in an implicit fashion using Lagrange interpolands $l^{(p)}_j$ of order $p$:
\begin{equation}\label{eq:interp}
u^N(t_n) \approx I^{(p)}u^N_n = \sum^p_{j=0} l^{(p)}_j u^N_{n+j-p}.
\end{equation}
Substituting (\ref{eq:interp}) into (\ref{eq:conv_short_long}), and evaluating the convolution integrals, we obtain the following approximations for the local and history parts, respectively, as
\begin{equation} \label{eq:discrlocal}
\mathcal{L}^{(n,\alpha)} u = \sum^p_{j=0} d^{(p)}_j u_{n+j-p},
\end{equation}
\begin{equation} \label{eq:discrhist}
\mathcal{H}^{(n,\alpha)} u = \sum^{n-p}_{j=0} \left(b^{(1)}_{n-1-j} u_j + b^{(2)}_{n-1-j} u_{j+1} + b^{(3)}_{n-1-j} u_{j+2} \right),
\end{equation}
where the corresponding $\alpha$- and $\Delta t$- dependent coefficients $d^{(p)}_j$ and $b^{(1)}_j,\, b^{(2)}_j,\,b^{(3)}_j$ are presented in \ref{Ap:Discr_Coefs}. In our computations, we make use of $p=1$ (piece-wise linear approximation) for the local part $\mathcal{L}^{(1,\alpha)} u$, when $n=1$ (first time-step), and $p=2$ for $n > 1$ (subsequent time-steps). The fractional RL derivative can be obtained from the above discretization by setting $-1 < \alpha < 0$. Therefore:
\begin{equation}\label{eq:discr_RL}
\prescript{RL}{0}{}\mathcal{D}^\alpha u^N(t)\big|_{t = t_n} \approx \mathcal{L}^{(n,-\alpha)} u^N + \mathcal{H}^{(n,-\alpha)} u^N.
\end{equation}
Finally, substituting (\ref{eq:discr_RL}) into (\ref{eq:discretized_initial}), and recalling (\ref{eq:discrlocal}), we obtain the discretized form for our FDE:
\begin{equation}\label{eq:discr_2}
	\sum^p_{j=0} d^{(p)}_j u^N_{n+j-p} + \mathcal{H}^{(n,-\alpha)} u^N + \sum^M_{j=1} W_{j,n}(\boldsymbol{\sigma}) \left(u^N_j - u_0 \right) = f^{data}_n,\quad n=1,\dots,\tilde{N}.
\end{equation}

\begin{remark}
Among a variety of available schemes, adopted the discretization of fractional operators introduced in \cite{Zheng2017} using Lubich's correction terms. However, in \cite{Zheng2017}, the authors do not perform a fully implicit computation of the second term on the left-hand-side of (\ref{eq:discretized_initial}) when $M>1$. Instead, they obtain the solutions for $n=1,\, \dots,\,M$ using a fine, uniform time-grid with time-step size $\tau = \Delta t^2$ and $M=1$ correction term. We observe that using such a fine initial time-grid might be computationally expensive and ineffective when strong singularities are present (\textit{e.g.} $u(t) = t^{\sigma^*}$, with $0 < \sigma^* < 1-\alpha$). In our developed scheme, we treat such term in a fully implicit fashion, where we solve a small linear system with $M$ unknowns to obtain $u^N_n$ for $n=1,\dots,\,M$. This ensures a proper inclusion of all $M$ singularities in all time-steps near $t=0$.
\end{remark}

We present here the developed finite-difference scheme to solve the discretized FDE (\ref{eq:discr_2}) for the particular case of $M=3$. Here, we solve a small linear system for the fully-implicit computation at the first 3 time-steps, that is, $t = \lbrace \Delta t,\,2\Delta t,\,3 \Delta t \rbrace$, which is obtained by expanding (\ref{eq:discr_2}) for each of the time-steps, where we use $p=1$ for the first one, with $\mathcal{H}^{(1,-\alpha)}u^N = 0$, and $p=2$ for the subsequent ones. The expansions are presented as follows:

\noindent
\textbf{First time-step} $t=\Delta t$:
\begin{equation}\label{eq:s1}
\underbrace{\left(d^{(1)}_1 + W_{1,1}\right)}_\text{$A_{11}$} u^N_1 +
\underbrace{\left(W_{2,1}\right)}_\textbf{$A_{12}$}u^N_2 + 
\underbrace{\left(W_{3,1}\right)}_\text{$A_{13}$} u^N_3 = f^{data}_1 - r_1 u_0,
\end{equation}
with,
\begin{equation}
r_1 = d^{(1)}_0 - \left(W_{1,1} + W_{2,1} + W_{3,1}\right).
\end{equation}

\noindent
\textbf{Second time-step} $t=2\Delta t$:
\begin{equation}\label{eq:s2}
\underbrace{\left(d^{(2)}_1 + b^{(2)}_0 + W_{1,2}\right)}_\text{$A_{21}$} u^N_1 + \underbrace{\left(d^{(2)}_2 + b^{(3)}_0 + W_{2,2}\right)}_\text{$A_{22}$} u^N_2 + \underbrace{\left(W_{3,2}\right)}_\text{$A_{23}$} u^N_3 = f^{data}_2 - r_2 u_0,
\end{equation}
with,
\begin{equation}
r_2 = d^{(2)}_0 + b^{(1)}_0 - \left(W_{1,2} + W_{2,2} + W_{3,2}\right).
\end{equation}

\noindent
\textbf{Third time-step} $t=3\Delta t$:
%
\begin{align}\label{eq:s3}
& \underbrace{\left(d^{(2)}_0 + b^{(1)}_0 + b^{(2)}_1 + W_{1,3}\right)}_\text{$A_{31}$} u^N_1 +
\underbrace{\left(d^{(2)}_1 + b^{(2)}_0 + b^{(3)}_1 + W_{2,3}\right)}_\text{$A_{32}$} u^N_2 \nonumber \\
& \qquad \qquad {} + \underbrace{\left(d^{(2)}_2 + b^{(3)}_0 + W_{3,3}\right)}_\text{$A_{33}$} u^N_3 = f^{data}_3 - r_3 u_0,
\end{align}
with,
\begin{equation}\label{eq:r3}
r_3 = b^{(1)}_1 - \left(W_{1,3} + W_{2,3} + W_{3,3}\right).
\end{equation}

Therefore, from (\ref{eq:s1})-(\ref{eq:r3}), we solve the following linear system for $n \le 3$:
\begingroup
\renewcommand*{\arraystretch}{1.25}
\begin{equation}\label{eq:u13}
\begin{bmatrix*}[l]
A_{11} & A_{12} & A_{13} \\
A_{21} & A_{22} & A_{23} \\
A_{31} & A_{32} & A_{33} 
\end{bmatrix*}
\begin{bmatrix*}
u^N_1 \\ u^N_2 \\ u^N_3
\end{bmatrix*} = 
\begin{bmatrix*}
f^{data}_1 \\ f^{data}_2 \\ f^{data}_3
\end{bmatrix*}
-
\begin{bmatrix*}
r_1 \\ r_2 \\ r_3
\end{bmatrix*} u_0, \quad \mathrm{if}\,\,n \le 3.
\end{equation}
\endgroup
\begin{remark}
To obtain the solutions for $M=2$, we only need to remove the third row/column of the coefficient matrix $\mathbf{A}$, and set $W_{1,3} = W_{3,1} = W_{2,3} = W_{3,2} = 0$ in (\ref{eq:s1})-(\ref{eq:r3}). Also, we observe that the correction terms lead to a full matrix of coefficients, which is reduced to a lower-triangular form for the uncorrected case. 
\end{remark}

For the remaining time-steps $(n > 3)$, $u^N_1,\,u^N_2,\,u^N_3$ are known and we solve for $u^N_n$ in the usual (but still implicit) time-stepping fashion, as follows:
\begin{equation}\label{eq:ug3}
	u^N_{n} = \left(f^{data}_n - \sum^1_{j=0} d^{(2)}_{j} u^N_{n-2+j} - \mathcal{H}^\alpha_n - \sum^M_{j=1} W_{j,n} (u^N_j - u_0) \right) / d^{(2)}_2, \,\, \mathrm{if}\,\, 3 < n \le \tilde{N}.
\end{equation}

\begin{remark}
The current formulation can be extended for a larger number of correction terms; however, using $M>9$ will incur in an ill-conditioned system for (\ref{eq:discretized_correction}). This fact was first analyzed by Diethelm .et.al \cite{Diethelm2006} and later numerically shown by Zeng .et.al \cite{Zeng2017CMAME}, which would incur in significant errors when computing the weights $W_{n,j}$ and consequently propagate the errors to the operator discretization. In that sense, we choose $M=3$ to capture the most critical singularities while keeping a small condition number and small residual for the linear system (\ref{eq:discretized_correction}). 
\end{remark}

\subsubsection{Stage-I Algorithm}
\label{Sec:Self-framework}

%
%

The Stage-I framework is described in Algorithm \ref{alg:1}, which learns about the singularities of the solution using small $\tilde{N}$. The $M$ captured singularities $\boldsymbol{\sigma}$ are utilized to initialize an FDE solver in Stage-II (see Section \ref{Sec:Stage-II}).

\renewcommand\thealgorithm{I-1}
\begin{algorithm}[t!]
	\caption{Stage-I Self-Singularity-Capturing Scheme for FDEs.}
	\label{alg:1}
	\begin{algorithmic}[1]
		\STATE{Known data about $u^{data}$, $f^{data}$ for $\tilde{N}$ time-steps of size $\Delta t$.}
		\STATE{Set $M = 1$, initial guess $\sigma^{(0)} = 0$ and numerical tolerances $\epsilon$, $\epsilon_1$.}
		\WHILE{$M \leq 3$}
		\STATE{Estimate $\sigma^{(k)}_i$, $i=1, \dots, M$ using Algorithm \ref{alg:1} and obtain the error $E(\boldsymbol{\sigma}^{(k)})$.}
		\IF{$E(\boldsymbol{\sigma}^{(k)}) < \epsilon$}
		\STATE{The singularities were captured within the specified tolerance $\epsilon$.}
		\STATE{\textbf{break}}
		\ELSE
		\STATE{Additional correction terms needed: Set $M = M + 1$.}
		\IF{$M == 2$} 
		\STATE{Set initial guess $\boldsymbol{\sigma}^{(0)} = [0,\,\,\, \sigma^{(k)}_1]^T$.}
		\ELSIF{$M == 3$} 
		\STATE{Set initial guess $\boldsymbol{\sigma}^{(0)} = [\sigma^{(k)}_1,\,\,\, \sigma^{(k)}_2,\,\,\, (\sigma^{(k)}_1 + \sigma^{(k)}_2)/2]^T$.}
		\ENDIF
		\ENDIF
		\ENDWHILE
		\RETURN{$M$, $\boldsymbol{\sigma}$, $E(\boldsymbol{\sigma})$.}
	\end{algorithmic}
\end{algorithm}

Algorithm \ref{alg:2} summarizes the main steps for the numerical scheme to capture the singularities of $u^{data},\,f^{data}$ for diminutive times at fixed $M$. 

\renewcommand\thealgorithm{I-2}
\begin{algorithm}[t!]
\caption{Stage-I Singularity Capturing Algorithm (with $M$ correction terms).}
\label{alg:2}
\begin{algorithmic}[1]
\STATE{Initial guess $\boldsymbol{\sigma}^0$, $\gamma^0$, compute $W_{n,j}$ using (\ref{eq:discretized_correction}) and $E(\boldsymbol{\sigma}^0)$ using (\ref{eq:Cost}), (\ref{eq:u13}), (\ref{eq:ug3}).}
\WHILE{$E(\boldsymbol{\sigma}^k) > \epsilon$ \AND $||\nabla E(\boldsymbol{\sigma}^k)|| > \epsilon_1$}
\STATE{Compute the perturbation $\boldsymbol{\sigma}^k + i\Delta \sigma \mathbf{e}_j$ for $j=1,\,2,\,\dots,\,M$.}
\STATE{Compute $\nabla E(\boldsymbol{\sigma}^k)$ using (\ref{eq:complex}), (\ref{eq:discretized_correction}), (\ref{eq:Cost}), (\ref{eq:u13}), (\ref{eq:ug3}).}
\STATE{Compute $\gamma^k$ using (\ref{eq:step}) and update $\boldsymbol{\sigma}^{k+1} = \boldsymbol{\sigma}^k - \gamma^k \nabla E(\boldsymbol{\sigma}^k)$.}
\ENDWHILE
\RETURN{$\boldsymbol{\sigma}$, $E(\boldsymbol{\sigma})$.}
\end{algorithmic}
\end{algorithm}

\subsubsection{Computational Complexity of Stage-I}
\label{Sec:Complexity}

We recall that in the presented scheme, the error function (\ref{eq:Cost}) is evaluated $N_{it}$ times, which is the number of iterations until convergence of the minimization scheme. The complexity of this error function is dominated by the time-integration of the numerical solution $u^N$ over $\tilde{N}$ time-steps. For the first $M$ time-steps, we use (\ref{eq:u13}) and solve the corresponding linear system, which costs $\mathcal{O}(M^3)$. For the remaining $\tilde{N}-M$ time-steps we use (\ref{eq:ug3}), where the dominant cost is due to the direct numerical evaluation of the fractional derivatives, that is, $\mathcal{O}((\tilde{N}-M)^2)$. Therefore, the computational complexity of the entire scheme is $\mathcal{O}( (M^3 + (\tilde{N}-M)^2) N_{it} )$. However, here we use $M \le 3$, and therefore the asymptotic complexity becomes $\mathcal{O}( \tilde{N}^2 N_{it})$. Furthermore, $\tilde{N}$ is assumed to be small due to the short-time data, and we show in Section \ref{Sec:Numerical} that we are able to capture singularities with $\tilde{N} = \mathcal{O}(M)$, which makes the presented scheme numerically efficient, as long as a large number of iterations for convergence is not required. 

\begin{remark}
	We observe that since $\tilde{N}$ is small, there is no need to use fast time-stepping schemes in Stage-I, since the break-even point between fast and direct schemes usually lies in a range of moderate number of time-steps (\textit{e.g.} about $\mathcal{O}(10^4)$ for the fast scheme in \cite{Zheng2017}). Therefore, for small $\tilde{N}$, fast schemes would decrease the performance of Stage-I, and would only be beneficial for Stage-II.
\end{remark}

\subsection{Stage-II: Integration of FDEs with the Captured Singularities}
\label{Sec:Stage-II}

Once the multi-singular behavior of the solution is captured through the framework presented in Section \ref{Sec:Self-framework} and Algorithm \ref{alg:1}, a variety of numerical methods for FDEs can be implemented (\textit{e.g.} other finite-difference and fast convolution schemes), incorporating the captured singularities. Here, to solve the FDE (\ref{eq:problemdef}), we use the developed implicit finite-difference scheme presented in Section \ref{Sec:FDun}, and given by expressions (\ref{eq:u13}) and (\ref{eq:ug3}) for $t \in \Omega = [t_0, t_N]$, using $N$ time-steps of size $\Delta t$ (not necessarily the same time-step size as Stage-I). The solution procedure is described in Algorithm \ref{alg:3}.


The computational complexity of Stage-II depends exclusively on the employed numerical discretization for the FDE. We utilize the same time-integration scheme as in Stage-I, for $N$ time-steps. Therefore, recalling (\ref{eq:u13}) and (\ref{eq:ug3}), we have a complexity of $\mathcal{O}( M^3 + (N-M)^2)$. However, since $M << N$, the asymptotic computational complexity of Stage-II becomes $\mathcal{O}(N^2)$.


\section{Numerical Results}
\label{Sec:Numerical}

We start with Stage-I with a computational analysis of the correction weights (see Section \ref{Sec:Num_Weights}), followed by the capturing scheme for particular case of a single time-step and correction term (see Section \ref{Sec:Num_simplified}). Then, we test the singularity-capturing Algorithm \ref{alg:2} for $M = 1,\,2$ and the self-capturing Algorithm \ref{alg:1} for $M=3$ (see Section \ref{Sec:Num_Mterms}). We also utilize the two-stage framework with random singularities and compare the accuracy of the entire time-integration framework with the captured singularities against \textit{ad-hoc} choices through a comparison between the obtained error functions $E(\boldsymbol{\sigma})$ (see Section \ref{Sec:random_1termFODE}). In all aforementioned tests, if not stated otherwise, we use the method of fabricated solutions, that is, we take $u^{data} = u^{ext}$ and $f^{data} = f^{ext}$, where we assume the following solution:
\begin{equation}\label{eq:udelta}
u^{ext}(t) = \sum^{S}_{j=1} t^{\sigma^*_j}, \quad \sigma^*_j \in \mathbb{R}^+,
\end{equation}
where $S$ denotes the number singularities/terms in $u^{ext} (t)$, and $\sigma^*_j$ represents prescribed singularities to be captured. From the defined analytical solution, we obtain the forcing term using a direct calculation, given by \cite{Podlubny99}:
\begin{equation}\label{eq:fdelta}
f^{ext}(t) = \sum^{S}_{j=1} \frac{\Gamma(1+\sigma^*_j)}{\Gamma(1+\sigma^*_j - \alpha)} t^{\sigma^*_j - \alpha}.
\end{equation}

Furthermore, we test Stage-I of the framework for a multi-term FDE (see Section \ref{Sec:multi-term}) and the two-stage framework for a one-term FDE with a \textit{singular-oscillatory} solution (see Section \ref{Sec:Oscillatory}).

\renewcommand\thealgorithm{II}
\begin{algorithm}[t!]
	\caption{Complete solution framework for FDEs.}
	\label{alg:3}
	\begin{algorithmic}[1]
		\STATE{\textbf{Stage-I: Self-Singularity-Capturing}}
		\STATE{Given data for $u^{data}$, $f^{data}$ for $\tilde{N}$ time-steps, capture $\boldsymbol{\sigma}$ using Algorithm \ref{alg:1}.}
		\STATE{\textbf{Stage-II: Numerical solution of the FDE}}
		\STATE{Let $N$ be the total number of time-steps. Compute $W_{n,j}$ using (\ref{eq:final_weights}).}
		\STATE{Solve $u^N_n$ for $n = 1,\dots,M$ using (\ref{eq:u13}).}
		\FOR{$n=4$ until $N$}
		\STATE{Solve for $u^N_n$ using (\ref{eq:ug3}).}
		\ENDFOR
	\end{algorithmic}
\end{algorithm}

\subsection{Numerical Behavior of Correction Weights}
\label{Sec:Num_Weights}

\begin{figure}[t!]
	\centering
	\begin{subfigure}[b]{0.4\textwidth}
		\includegraphics[width=\columnwidth]{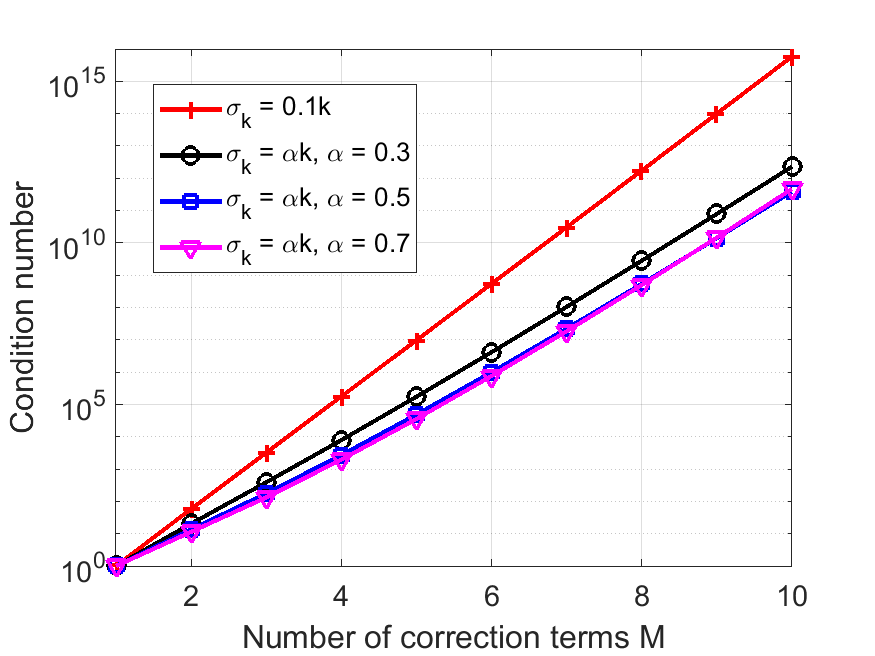}
		\caption{Condition number of $V$ \textit{vs} $M$. \label{fig:Cond_a}}	
	\end{subfigure}%
	~ 
	\begin{subfigure}[b]{0.4\textwidth}
		\includegraphics[width=\columnwidth]{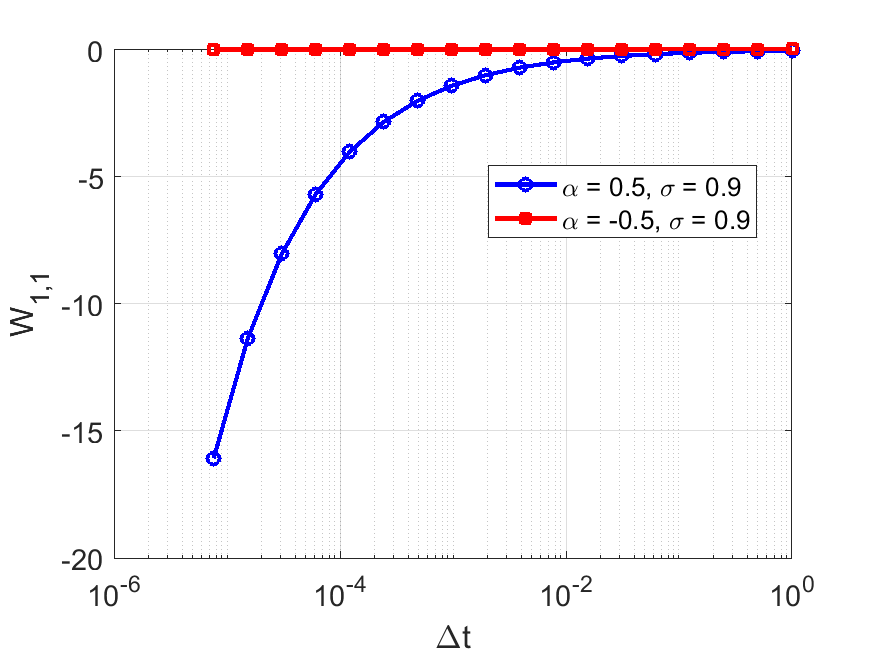}
		\caption{First correction weight $W_{1,1}$ \textit{vs} $\Delta t$. \label{fig:W11}}	
	\end{subfigure}
	\caption{\textit{(a)} Condition numbers for different choices of powers $\sigma_k$. We observe that the choice $\sigma_k = 0.1 k$ for unknown regularity of $u(t)$ \cite{Zheng2017} leads to the highest condition number, where at most $M=9$ correction terms can be used with double precision arithmetic. \textit{(b)} First correction weight $W_{1,1}$ with respect to $\Delta t$. \label{fig:Cond}}
\end{figure}

\begin{figure}[t!]
	\centering
	\begin{subfigure}[b]{0.4\textwidth}
		\includegraphics[width=\columnwidth]{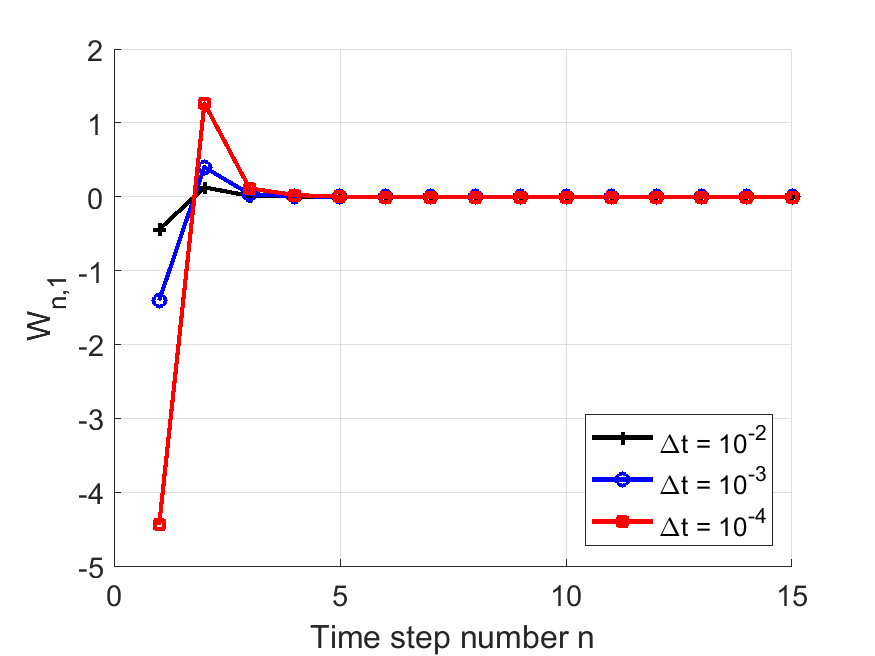}
		\caption{$\alpha = 0.5,\, \sigma_k = 0.9$. \label{fig:Walpha}}	
	\end{subfigure}%
	~
	\begin{subfigure}[b]{0.4\textwidth}
		\includegraphics[width=\columnwidth]{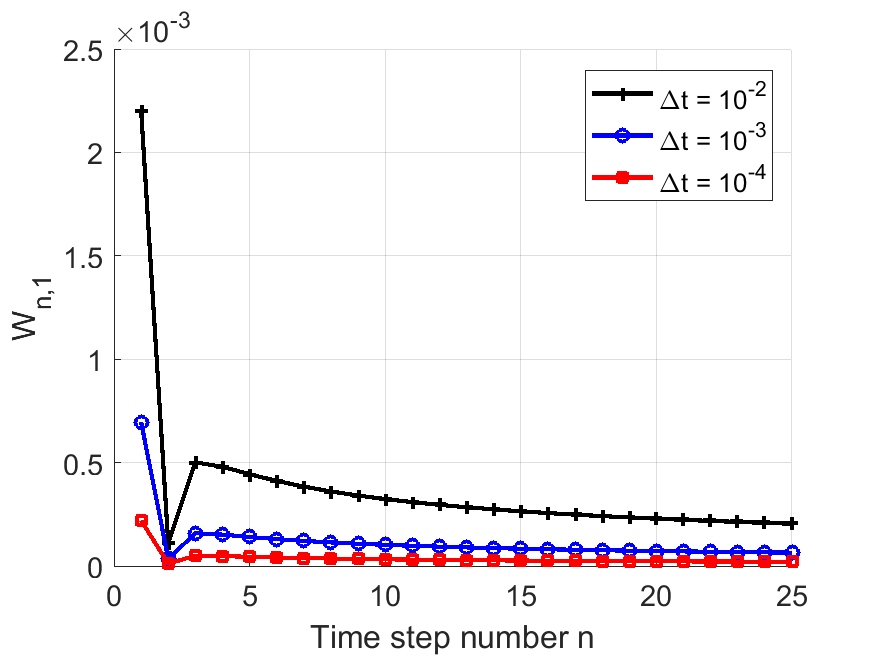}
		\caption{$\alpha = -0.5,\, \sigma_k = 0.9$. \label{fig:Wmalpha}}	
	\end{subfigure}
	\caption{Initial weights $W_{1,n}$ for different time-step sizes $\Delta t$, with corrections for fractional differentiation/integration \textit{(a)} $\alpha = 0.5$ (fractional differentiation), where the initial weights have a larger magnitude that quickly decreases in the first time-steps, indicating that they are most relevant near $t = 0$. \textit{(b)} $\alpha = -0.5$ (fractional integral). \label{fig:W_alpha}}
\end{figure}

We investigate the behavior of the initial correction weights $W_{n,j}$ and the condition number for the Vandermonde matrix $V$ presented in (\ref{eq:final_weights}). Figure \ref{fig:Cond_a} illustrates the condition numbers for $V_{M\times M}$ using different choices for $\sigma_k$, namely $\sigma = \alpha k$ when the regularity of $u(t)$ is known, and $\sigma_k = 0.1 k$ when it is unknown. We also illustrate the behavior of the correction weights using $M=1$ with respect to the time-step size $\Delta t$ and fractional order $\alpha$ (corresponding to fractional differentiation/integration) in Figures \ref{fig:W11} and \ref{fig:W_alpha}. We observe that for $\alpha = 0.5$ (fractional differentiation), $W_{1,1}$ increases in magnitude as $\Delta t$ decreases. On the other hand, for $\alpha = 0.5$ (fractional integration), it starts with $W_{1,1} \approx 0.2$ for $\Delta t = 1$, and decreases with $\Delta t$. For fractional integration, only positive correction weights are observed, that decrease slower with respect to $n$, but with all values decreasing to zero as $\Delta t$ decreases.

\subsection{Single Time-Step and Correction Term}
\label{Sec:Num_simplified}

\begin{figure}[t!]
	\centering
	\begin{subfigure}[b]{0.4\textwidth}
		\includegraphics[width=\columnwidth]{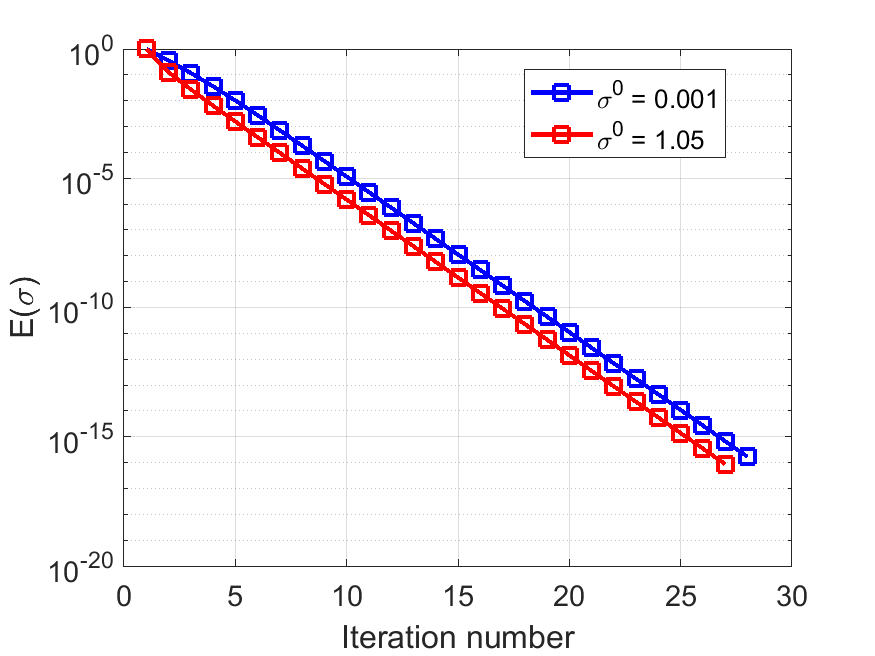}
		\caption{Error \textit{vs} iteration number.}		
	\end{subfigure}%
	~ 
	\begin{subfigure}[b]{0.4\textwidth}
		\includegraphics[width=\columnwidth]{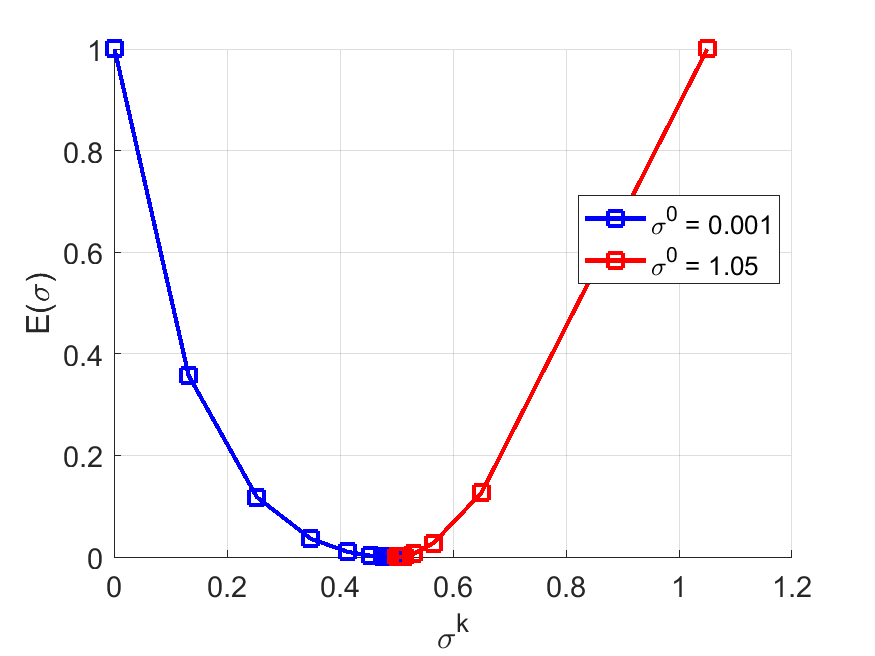}
		\caption{Error \textit{vs} iterated singularity}		
	\end{subfigure}
	\caption{Convergence behavior with $\sigma^* = 0.5$ and initial guesses $\sigma^0 = \lbrace 0.0001, 1.05 \rbrace$ and time-step size $\Delta t = 0.01$. \label{fig:f1}}
\end{figure}

\begin{figure}[t!]
	\centering
	\begin{subfigure}[b]{0.4\textwidth}
		\includegraphics[width=\columnwidth]{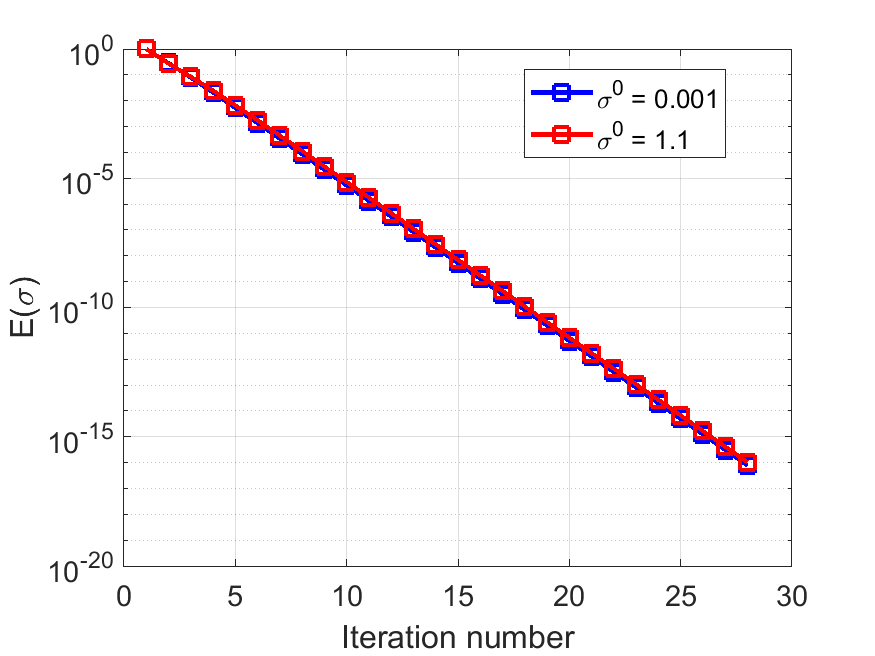}
		\caption{Error \textit{vs} iteration number}		
	\end{subfigure}%
	~ 
	\begin{subfigure}[b]{0.4\textwidth}
		\includegraphics[width=\columnwidth]{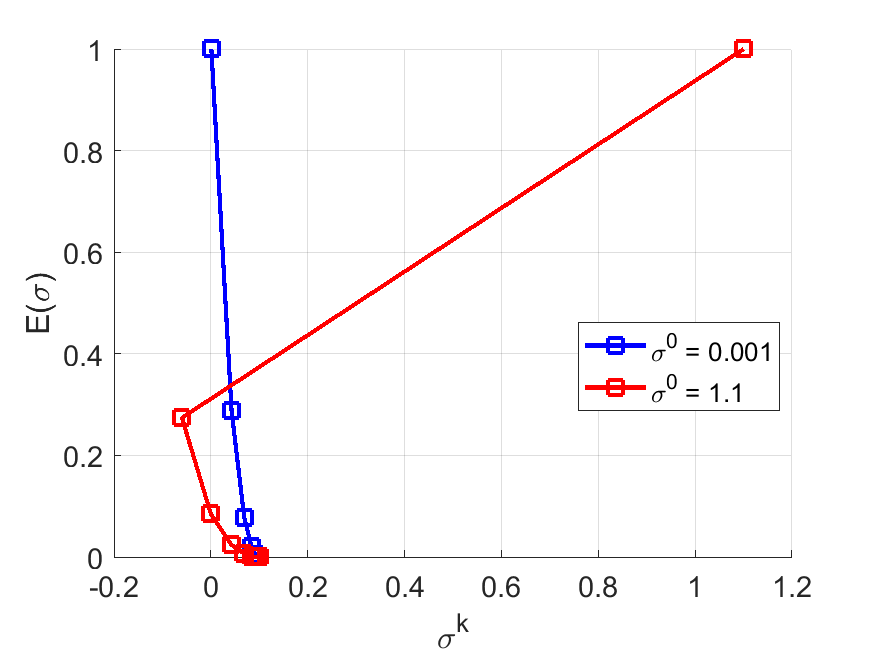}
		\caption{Error \textit{vs} iterated singularity}		
	\end{subfigure}
	\caption{Convergence behavior with $\sigma^* = 0.1$ and initial guesses $\sigma^0 = \lbrace 0.0001, 1.1 \rbrace$ and time-step size $\Delta t = 0.01$. \label{fig:f2}}
\end{figure}

\begin{figure}[t!]
	\centering
	\begin{subfigure}[b]{0.4\textwidth}
		\includegraphics[width=\columnwidth]{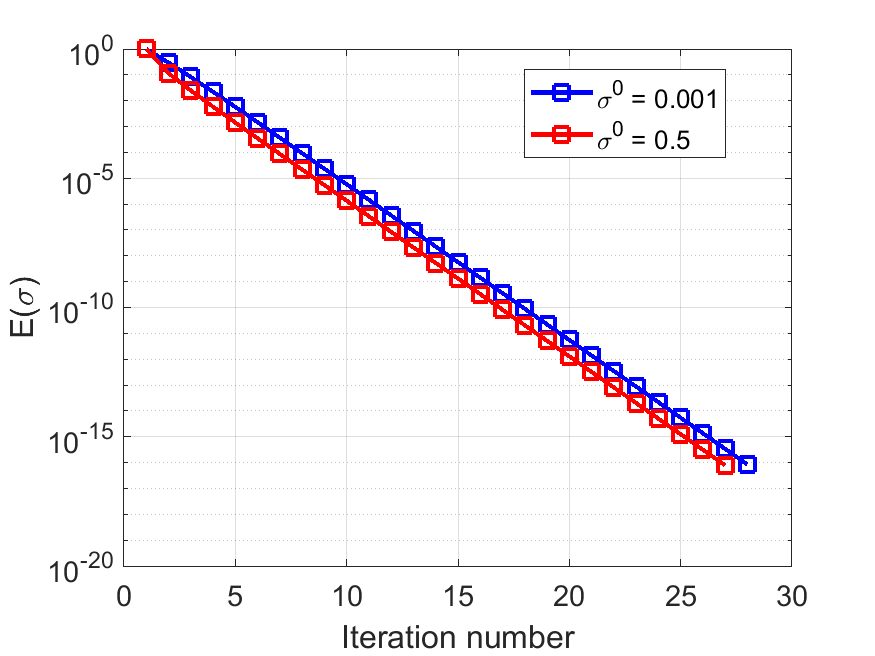}
		\caption{Error \textit{vs} iteration number}				
	\end{subfigure}%
	~ 
	\begin{subfigure}[b]{0.4\textwidth}
		\includegraphics[width=\columnwidth]{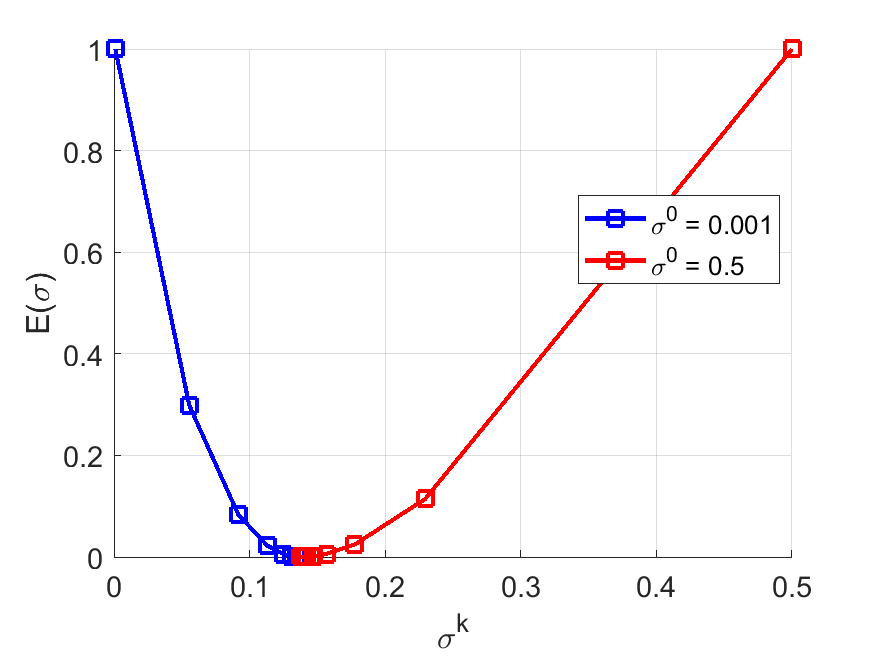}
		\caption{Error \textit{vs} iterated singularity}
	\end{subfigure}
	\caption{Convergence behavior with $\sigma^*_1 = 0.1$, $\sigma^*_2 = 0.2$, initial guesses $\sigma^0 = \lbrace 0.001, 0.5 \rbrace$ and time-step size $\Delta t = 0.01$. The converged value for the singularity $\sigma \approx 0.1377$ is between $\sigma^*_1, \sigma^*_2$. \label{fig:f3}}
\end{figure}

\begin{figure}[t!]
	\centering
	\begin{subfigure}[b]{0.4\textwidth}
		\includegraphics[width=\columnwidth]{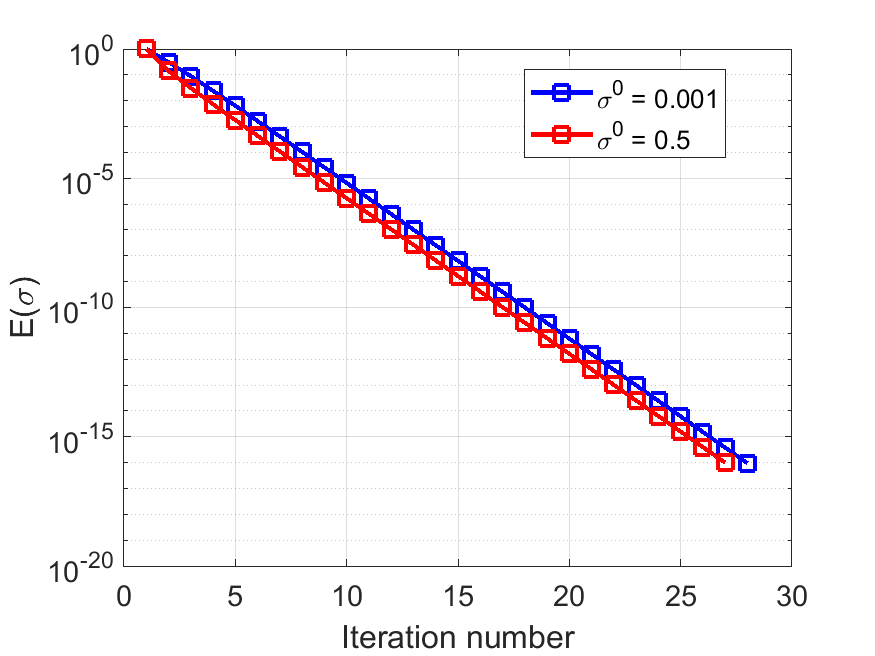}
		\caption{Error \textit{vs} iteration number}				
	\end{subfigure}%
	~ 
	\begin{subfigure}[b]{0.4\textwidth}
		\includegraphics[width=\columnwidth]{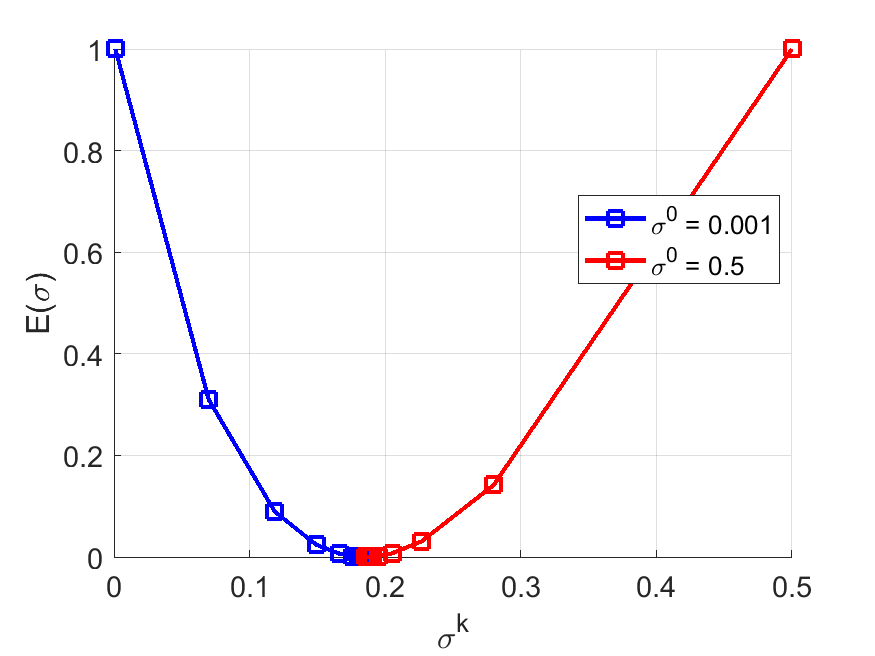}
		\caption{Error \textit{vs} iterated singularity}
	\end{subfigure}
	\caption{Convergence behavior with $\sigma^*_1 = 0.1$, $\sigma^*_2 = 0.3$, $\sigma^*_3 = 0.5$, initial guesses $\sigma^0 = \lbrace 0.001, 0.5 \rbrace$ and time-step size $\Delta t = 0.01$. The converged value for the singularity is $\sigma \approx 0.1856 $.\label{fig:f4}}
\end{figure}

We present the numerical results for a particular case of Stage-I in \ref{Ap:Single} for $S = 1$ singularity and $M = 1$ correction term, with a closed-form for the correction weights (see (\ref{eq:solution_W})). Let the time domain $\tilde{\Omega} = [0, \Delta t]$, with fractional-order and time-step size kept constant, respectively, at $\alpha = 0.5$ and $\Delta t = 0.01$ in this section. Figure \ref{fig:f1} presents the convergence behavior for $\sigma^*_1 = 0.5$, and Figure \ref{fig:f2} illustrates the obtained results using $\sigma^*_1 = 0.1$. We observe that in both cases, the scheme captures the power $\sigma^*$ from the analytical solution, with a relatively small number of iterations. Furthermore, we observe an overshooting in the iterative procedure for initial guess $\sigma^{(0)} = 1.1$, but nevertheless, the scheme still converges within machine precision. Figures \ref{fig:f3} and \ref{fig:f4} show, respectively, the results for $S=2$ with $\sigma^*_1 = 0.1$, $\sigma^*_2 = 0.2$ and $S=3$ with $\sigma^*_1 = 0.1$, $\sigma^*_2 = 0.3$, $\sigma^*_3 = 0.5$. We observe that the converged values for $\sigma$ are closer to $\sigma^*_1$, which is the most critical singularity for the chosen $u^{ext} (t)$.

\begin{figure}[t!]
	\centering
	\includegraphics[width=0.4\textwidth]{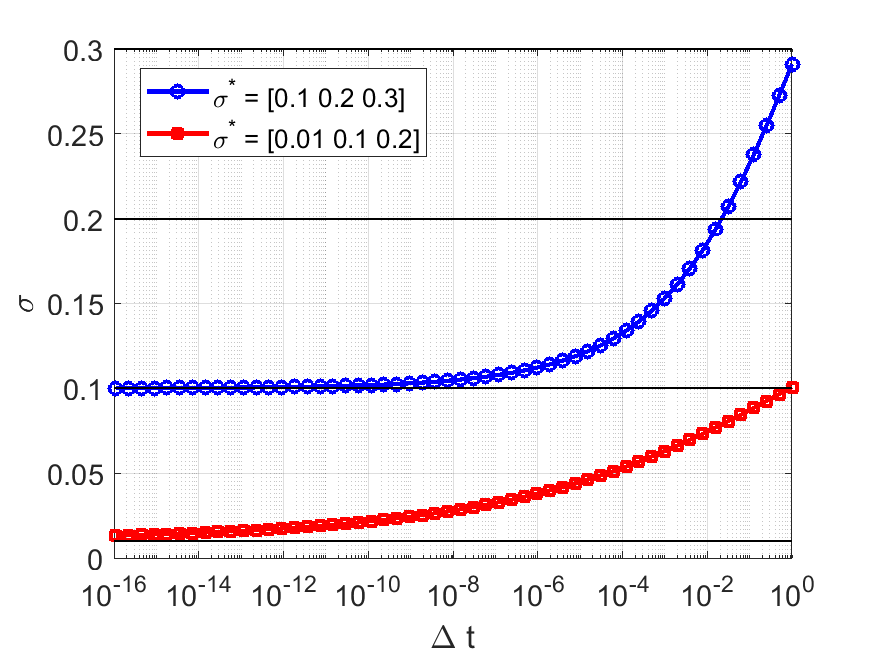}
	\caption{Converged values for the singularity $\sigma$ \textit{vs} time-step size $\Delta t$ with two choices of $\boldsymbol{\sigma}^* = \left[\sigma^*_1, \,\,\, \sigma^*_2, \,\,\, \sigma^*_3 \right]^T$. \label{fig:sigma_dt}}
\end{figure}

Since the converged values for $\sigma$ captured intermediate values between the defined singularities $\sigma^*_1, \sigma^*_2, \sigma^*_3$ for the choice of $\Delta t = 0.01$, we analyze the convergence behavior of $\sigma$ with respect to $\Delta t$. We present the obtained results in Figure \ref{fig:sigma_dt} for two different sets of singular values, and we observe for the defined range of $\Delta t$, that $\sigma$ lies between the singularities $\sigma^*_1$, $\sigma^*_2$, $\sigma^*_3$, and converge to the strongest singularity (in this case $\sigma^*_1$) as $\Delta t \to 0$.

\subsection{Singularity Capturing for $M$ Correction Terms}
\label{Sec:Num_Mterms}

We systematically test the capturing Algorithm \ref{alg:2} for multiple correction terms $M=1,2,3$ and multiple defined singularities $S = 1,2,3$ in $u^{ext}$. The tests are presented in an incremental fashion for $M$, where we show the capabilities of Algorithm \ref{alg:2} to capturing/approximating the singularities. We then demonstrate how the self-capturing framework defined in Algorithm \ref{alg:1} successfully determines the singularities for $S = 3$ and $M = 3$, where we compare the obtained errors with \textit{ad-hoc} choices for $\boldsymbol{\sigma}$ using random strong singularities. Throughout this section, we keep the fractional order fixed at $\alpha = 0.5$, as well as short time-domain $\tilde{\Omega} = [0, 1]$, and perturbation $\Delta \sigma = 10^{-14}$ for the complex step differentiation. For all cases, we define a tolerance $\epsilon = 10^{-15}$ for $E(\boldsymbol{\sigma})$ and unless stated otherwise, $\epsilon_1 = 10^{-14}$ for $||\nabla E(\boldsymbol{\sigma})||$. We choose a smaller tolerance for the norm of the error gradient (since the $E(\boldsymbol{\sigma})$ is defined with the norm $||.||^2_2$) to make sure that $E(\boldsymbol{\sigma})$ is always minimized before we introduce additional correction terms in the self-capturing Stage-I. We also use $\gamma^{0} = 10^{-3}$ for initialization in Algorithm \ref{alg:2}. For cases where $M = S$, we compute the component-wise relative error of the converged $\boldsymbol{\sigma}^k$, which we define as:
\begin{equation}
E^\sigma_j = \frac{\vert \sigma^*_j - \sigma^{(k)}_j \vert}{\vert \sigma^*_j \vert},\quad j=1,\dots, M,
\end{equation}
where $k$ denotes the iteration number when convergence is achieved, \textit{i.e.}, $E(\boldsymbol{\sigma}^{(k)}) < \epsilon$.

\subsubsection{One Correction Term}

\begin{figure}[t!]
	\centering
	\begin{subfigure}[b]{0.4\textwidth}
	\includegraphics[width=\columnwidth]{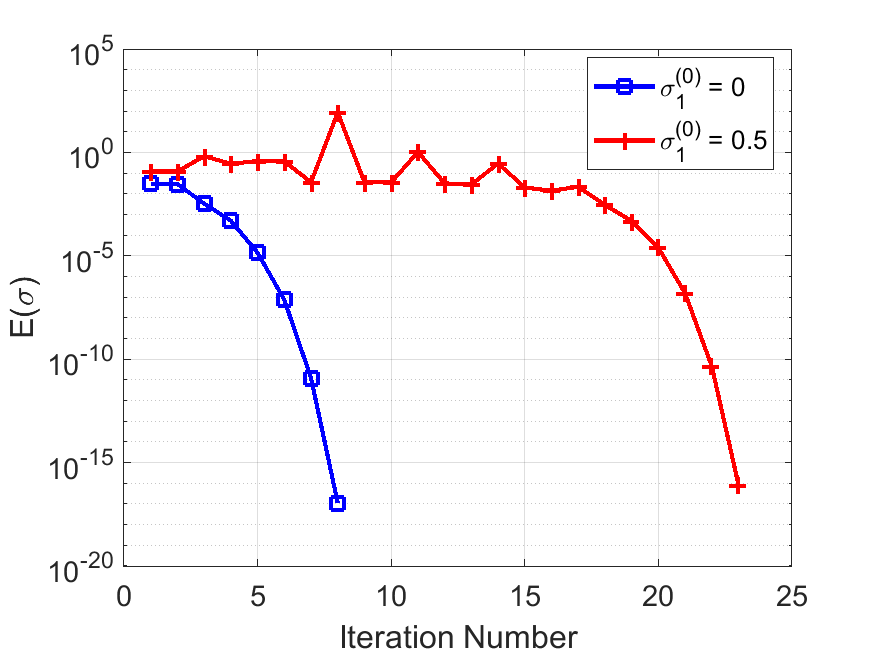}
		\caption{Error \textit{vs} iteration number}		
	\end{subfigure}%
	~ 
	\begin{subfigure}[b]{0.4\textwidth}
		\includegraphics[width=\columnwidth]{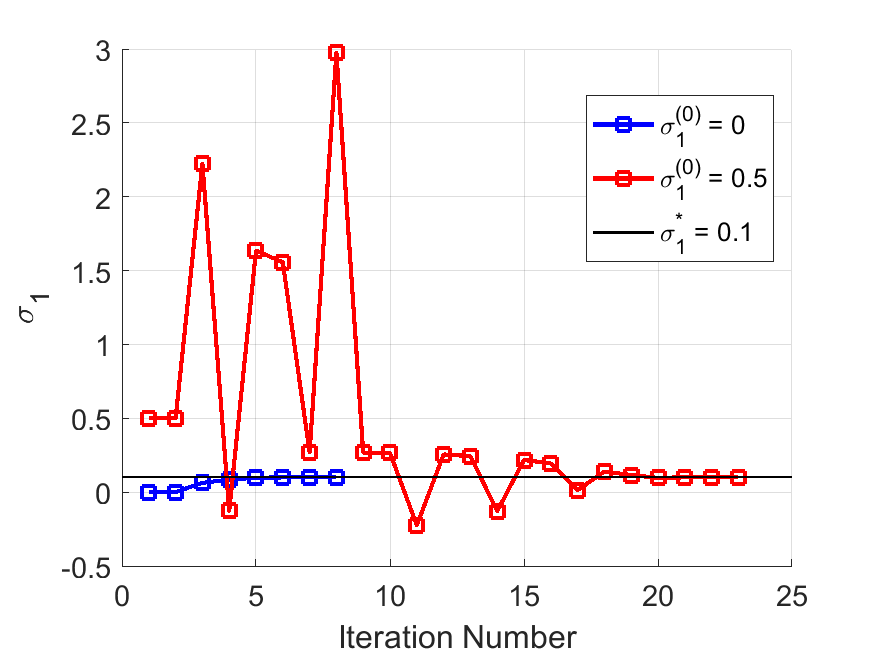}
		\caption{$\sigma$ \textit{vs} iteration number}		
	\end{subfigure}%
	\caption{Second-order asymptotic convergence behavior for $M=1$ correction term, $S = 1$ singularity with value $\sigma^*_1 = 0.1$. The obtained relative error for the converged $\sigma_1$ is $E^\sigma_1 = 2.290 \times 10^{-8}$.\label{fig:f5}}
\end{figure}

\begin{figure}[t!]
	\centering
	\begin{subfigure}[b]{0.4\textwidth}
	\includegraphics[width=\columnwidth]{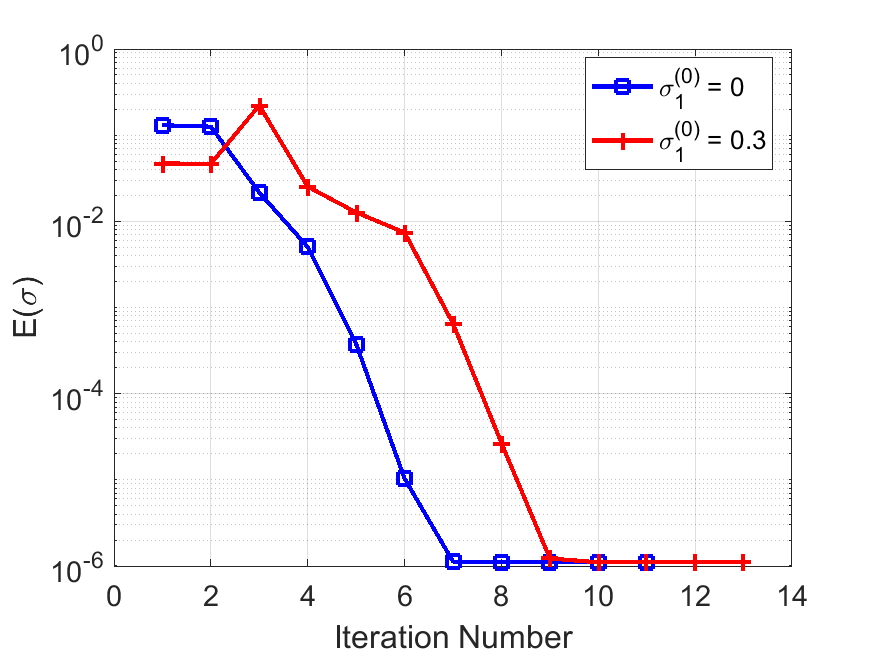}
		\caption{Error \textit{vs} iteration number.}		
	\end{subfigure}%
	~ 
	\begin{subfigure}[b]{0.4\textwidth}
		\includegraphics[width=\columnwidth]{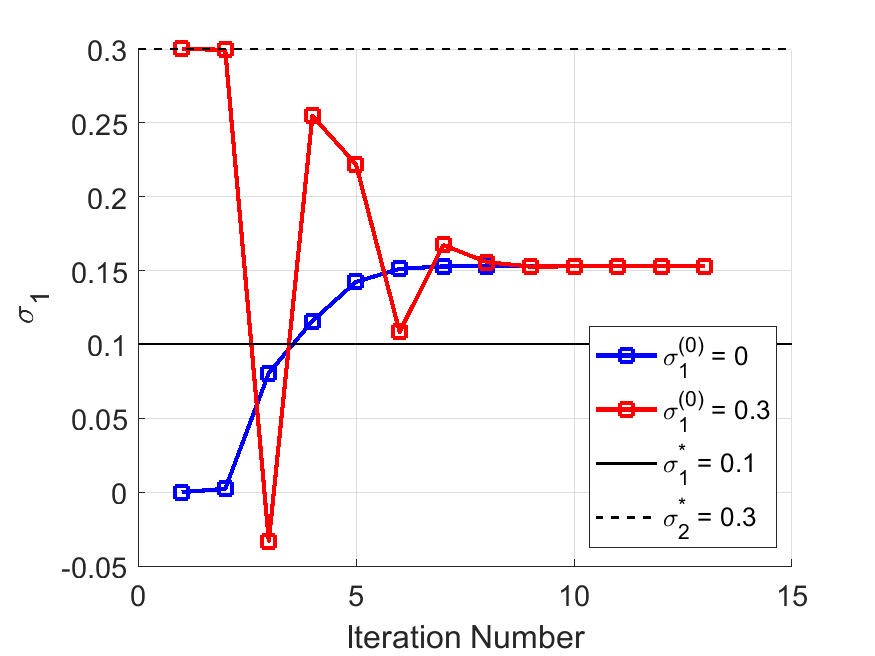}
		\caption{$\sigma$ \textit{vs} iteration number.}		
	\end{subfigure}%
	\caption{Convergence behavior for $M=1$, $S = 2$ with $\boldsymbol{\sigma}^* = \left[0.1,\,\,0.3\right]^T$. With the plateau in error and $\sigma$ due to gradient convergence within the tolerance $\epsilon_1$. The scheme converges to an ``\textit{intermediate}" value $\sigma \approx 0.153$, between the true singularities. \label{fig:f6}}
\end{figure}

\begin{figure}[t!]
	\centering
	\begin{subfigure}[b]{0.4\textwidth}
	\includegraphics[width=\columnwidth]{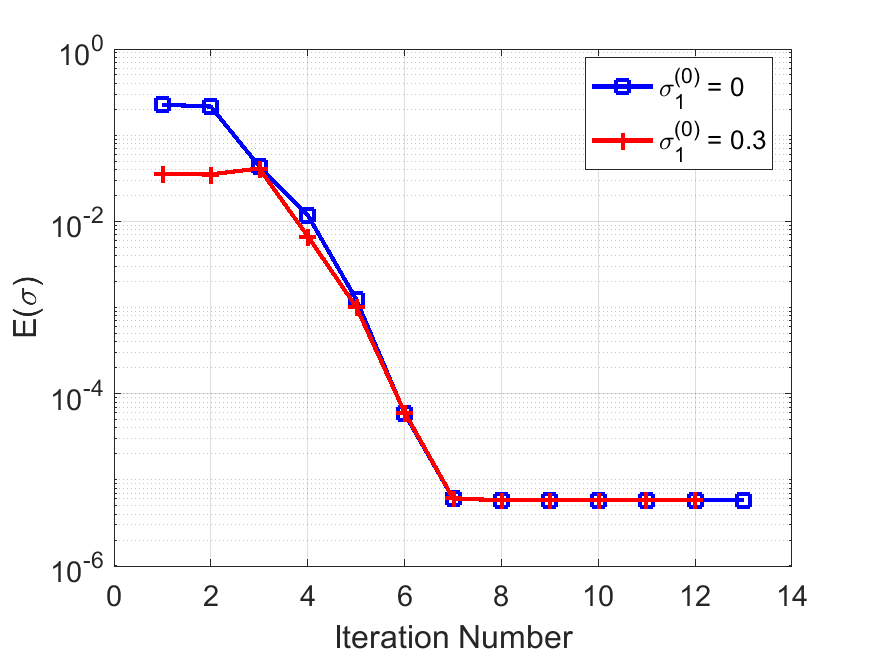}
		\caption{$\sigma$ \textit{vs} iteration number.}		
	\end{subfigure}%
	~ 
	\begin{subfigure}[b]{0.4\textwidth}
		\includegraphics[width=\columnwidth]{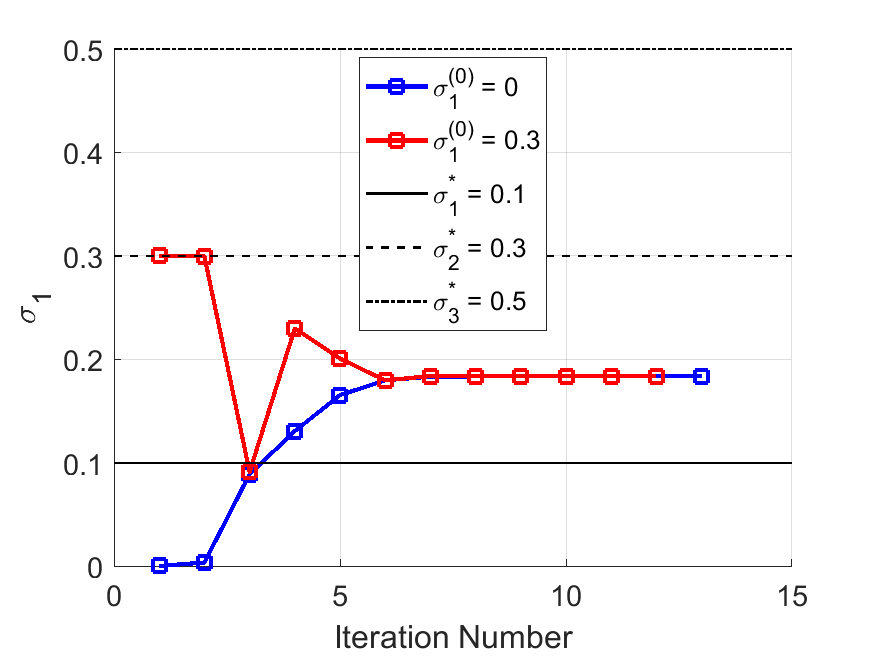}
		\caption{$\sigma$ \textit{vs} iteration number.}		
	\end{subfigure}%
	\caption{Convergence behavior for $M=1$, $S = 3$ with $\boldsymbol{\sigma}^* = \left[0.1,\,\, 0.3,\,\, 0.5\right]^T$. In for both initial guesses, the scheme converges to an ``\textit{intermediate}" value $\sigma \approx 0.184$, which is closer to the strongest true singularity $\sigma^*_1 = 0.1$.\label{fig:f7}}
\end{figure}

We consider $M = 1$ and $\tilde{N} = 100$ initial data-points. Figures \ref{fig:f5}, \ref{fig:f6} and \ref{fig:f7} present the obtained results, respectively, for $S = 1,\, 2,\,3$ singularities, using Algorithm \ref{alg:2}. We observe that, when $S=M$, we can capture singularities within machine precision. When $M<S$, we are still able to find a minimum for $E(\sigma)$ with an intermediate value for $\sigma_1$, similar to the results obtained for the simplified case in Section \ref{Sec:Num_simplified}. Furthermore, when $S,\,M \ge 2$, or the initial guess $\sigma^{(0)}$ is far from the true singularity, the iterative procedure assumes the typical zig-zag behavior of the gradient descent scheme.


\subsubsection{Two Correction Terms}

\begin{figure}[t!]
	\centering
	\begin{subfigure}[b]{0.4\textwidth}
	\includegraphics[width=\columnwidth]{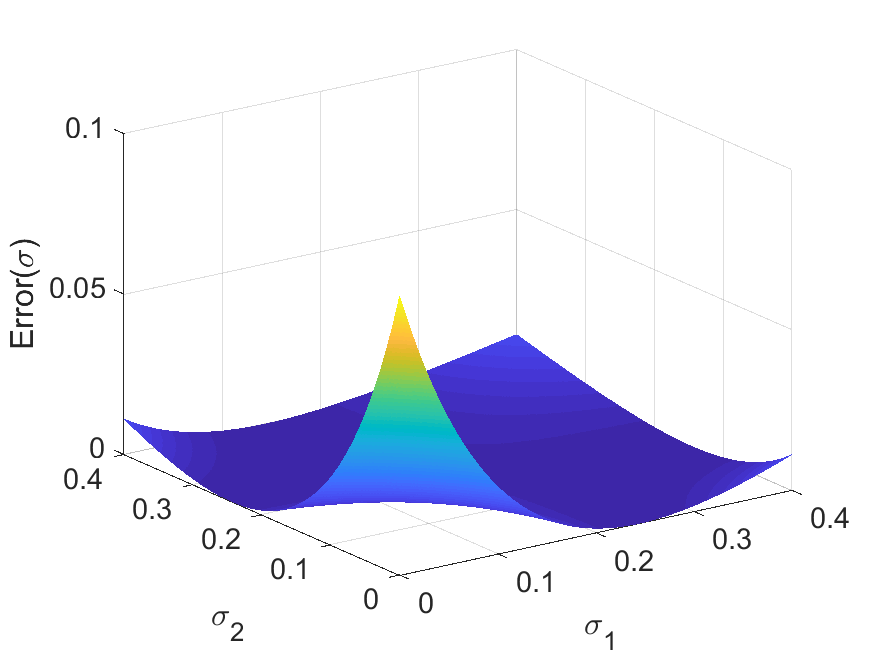}
		\caption{Linear error-axis.}		
	\end{subfigure}%
	~ 
	\begin{subfigure}[b]{0.4\textwidth}
		\includegraphics[width=\columnwidth]{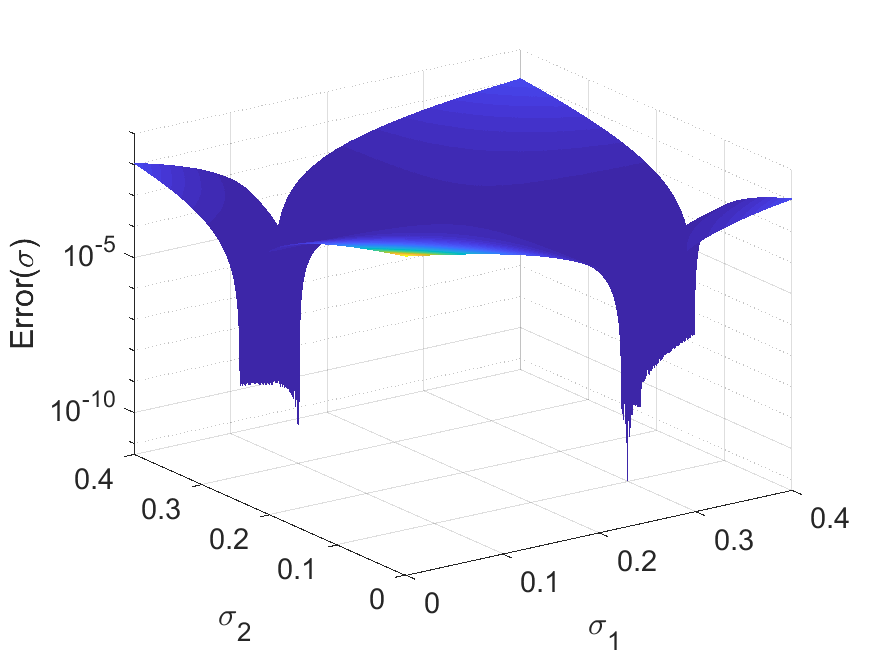}
		\caption{Logarithmic error-axis.}		
	\end{subfigure}%
	\caption{Error function $E(\boldsymbol{\sigma})$ for $M=2$, $S = 2$, with $\sigma^*_1 = 0.1$, $\sigma^*_2 = 0.3$. Two minima occur at $\boldsymbol{\sigma}^* = \left[0.1,\,\, 0.3\right]^T = \left[0.3,\,\, 0.1\right]^T$, corresponding to the true singularities. We set $\sigma_1 \neq \sigma_2$ to avoid a singular matrix in linear system (\ref{eq:discretized_correction}). \label{fig:f8}}
\end{figure}

\begin{figure}[t!]
	\centering
	\begin{subfigure}[b]{0.4\textwidth}
	\includegraphics[width=\columnwidth]{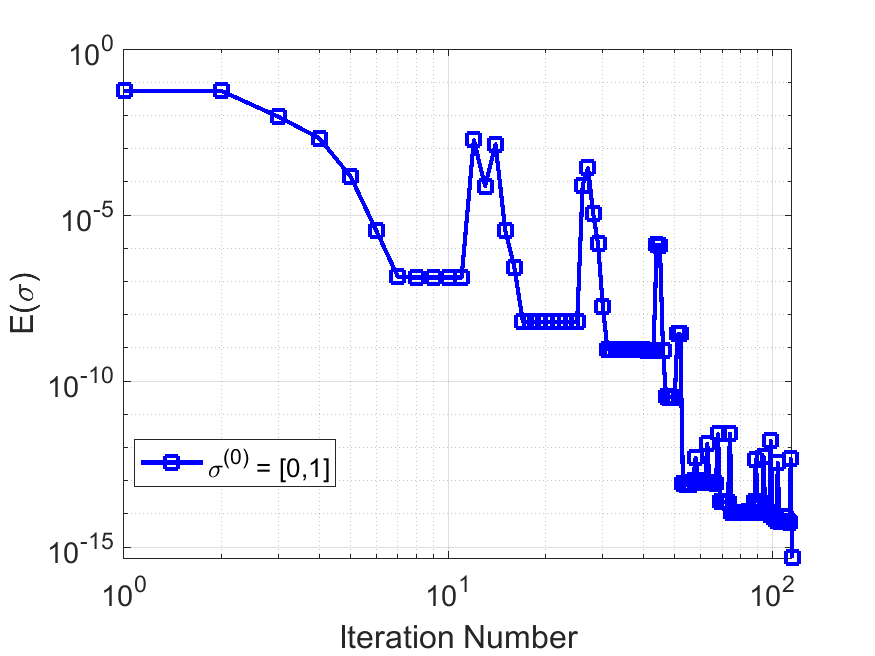}
		\caption{Error \textit{vs} iteration number.}		
	\end{subfigure}%
	~ 
	\begin{subfigure}[b]{0.4\textwidth}
		\includegraphics[width=\columnwidth]{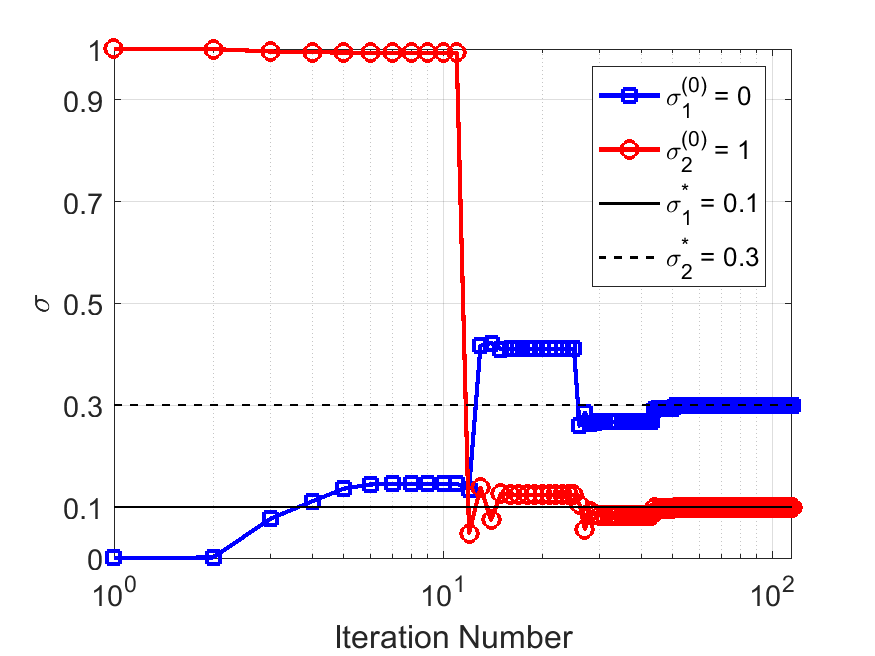}
		\caption{$\boldsymbol{\sigma}$ \textit{vs} iteration number.}		
	\end{subfigure}%
	\caption{Two correction terms $M=2$ and two singularities $S = 2$, with $\boldsymbol{\sigma}^* = \left[0.1,\,\, 0.3\right]^T$. The obtained converged values are $\boldsymbol{\sigma} = \left[0.299974,\,\, 0.099990\right]^T$, which correspond to component-wise relative errors $E^\sigma = \left[8.61,\,\, 9.94\right]^T \times 10^{-5}$. \label{fig:f9}}
\end{figure}

\begin{figure}[t!]
	\centering
	\begin{subfigure}[b]{0.4\textwidth}
	\includegraphics[width=\columnwidth]{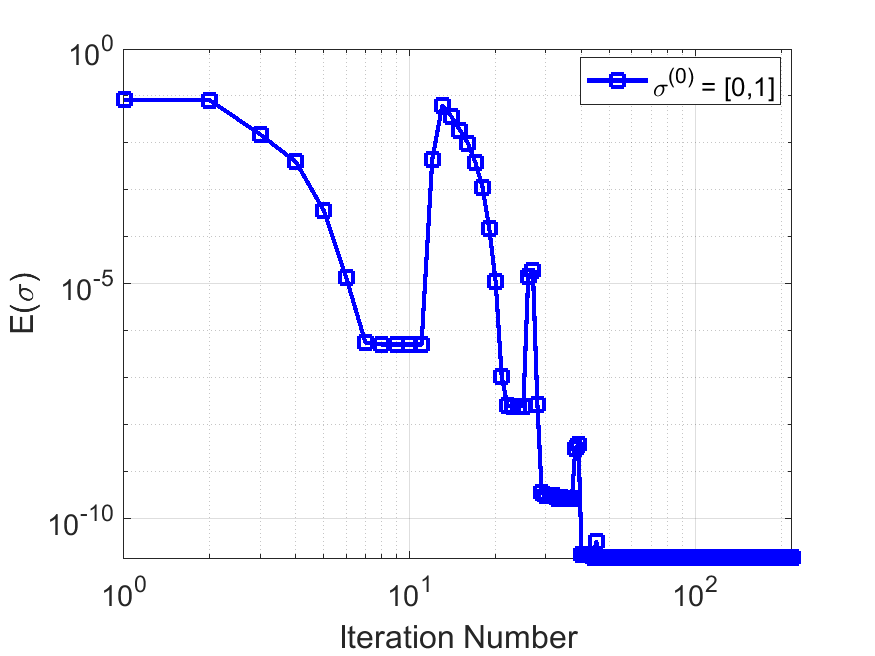}
		\caption{Error \textit{vs} iteration number.}		
	\end{subfigure}%
	~ 
	\begin{subfigure}[b]{0.4\textwidth}
		\includegraphics[width=\columnwidth]{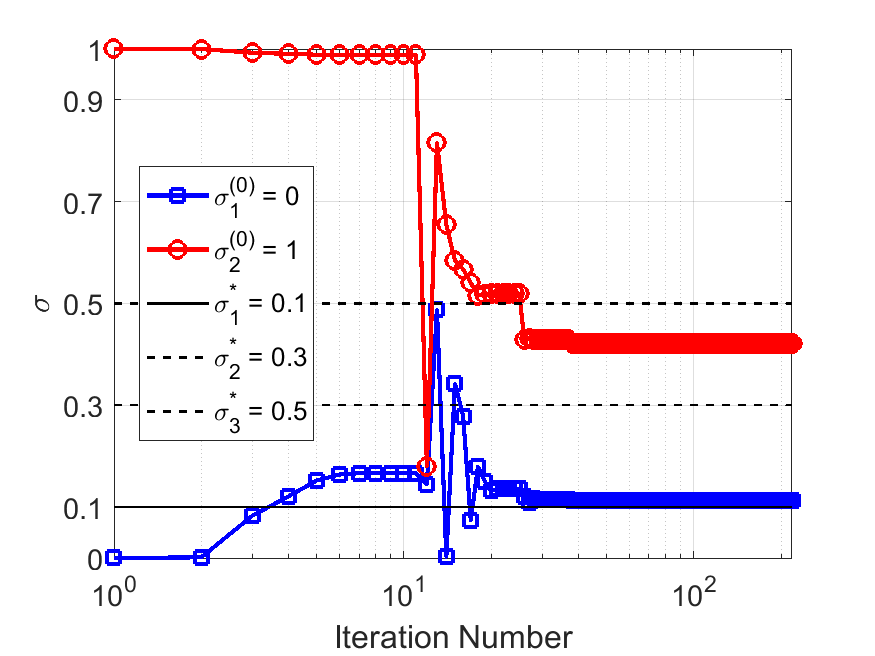}
		\caption{$\boldsymbol{\sigma}$ \textit{vs} iteration number.}		
	\end{subfigure}%
	\caption{Two correction terms $M=2$ and three singularities $S = 3$, with $\boldsymbol{\sigma}^* = \left[0.1,\,\, 0.3\right]^T$ and $\epsilon_1 = 10^{-11}$. The converged values are $\boldsymbol{\sigma} = \left[0.112635,\,\, 0.420495\right]^T$, which are closer to the stronger singularities $\sigma_1 = 0.1$ and $\sigma_3 = 0.5$. The approximate singularities still lead to an approximation error of $1.47 \times 10^{-11}$. \label{fig:f10}}
\end{figure}

We consider $M=2$, where we first analyze the behavior of the error function $E(\boldsymbol{\sigma})$ using $S = 2$ singularities (\textit{see Figure \ref{fig:f8}}). Then, we use Algorithm \ref{alg:2} to capture the singularities $\boldsymbol{\sigma}^* = \left[0.1,\,\,0.3\right]^T$ (\textit{see Figure \ref{fig:f9}}), where we observe a more pronounced zig-zag behavior for convergence, compared to $M=1$. We also obtain estimates of $\sigma_1,\,\sigma_2$ when dealing with $S = 3$ singularities, namely $\boldsymbol{\sigma}^* = \left[0.1,\,\,0.3,\,\,0.5\right]^T$ (\textit{see Figure \ref{fig:f10}}). Of particular interest, we observe that the obtained singularities lie in the range of the true singularities, which are useful observations to define initial guesses for cases with three correction terms and singularities.

\subsubsection{Three Correction Terms}

\begin{figure}[t!]
	\centering
	\begin{subfigure}[b]{0.4\textwidth}
	\includegraphics[width=\columnwidth]{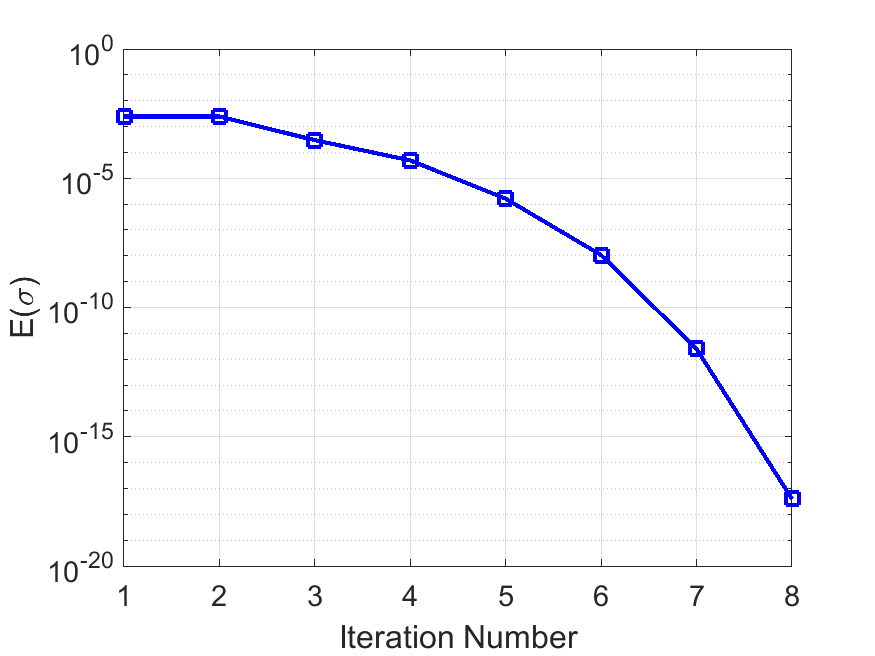}
		\caption{Error \textit{vs} iteration number.}		
	\end{subfigure}%
	~ 
	\begin{subfigure}[b]{0.4\textwidth}
		\includegraphics[width=\columnwidth]{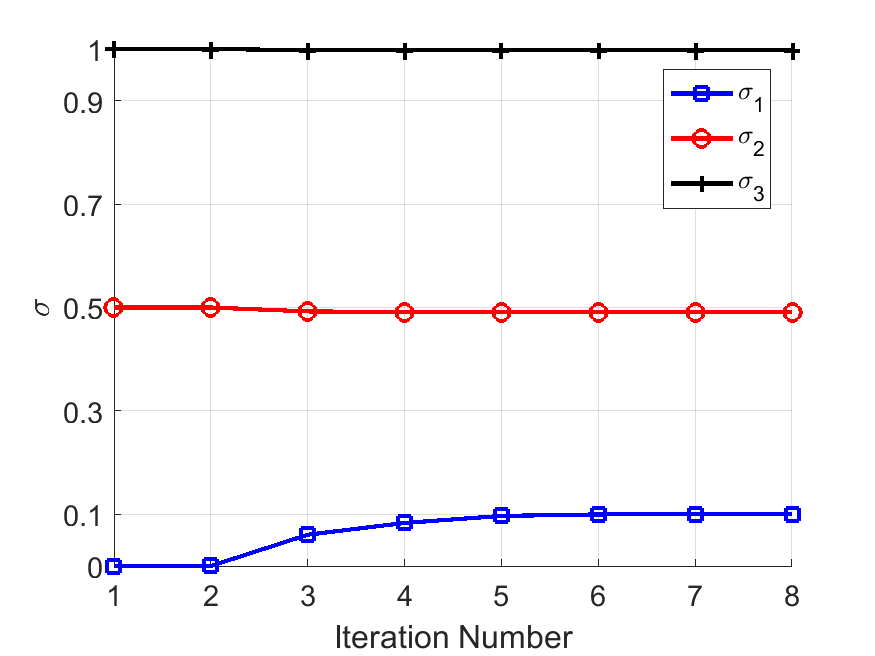}
		\caption{$\boldsymbol{\sigma}$ \textit{vs} iteration number}		
	\end{subfigure}%
	\caption{Results for $M=3$, $S = 1$. The obtained relative error for $\sigma_1$ is $E^\sigma_1 = 4.869\times 10^{-8}$. The components $\sigma_1,\,\sigma_2$ remained practically constant.\label{fig:f11}}
\end{figure}

\begin{figure}[t!]
	\centering
	\begin{subfigure}[b]{0.4\textwidth}
	\includegraphics[width=\columnwidth]{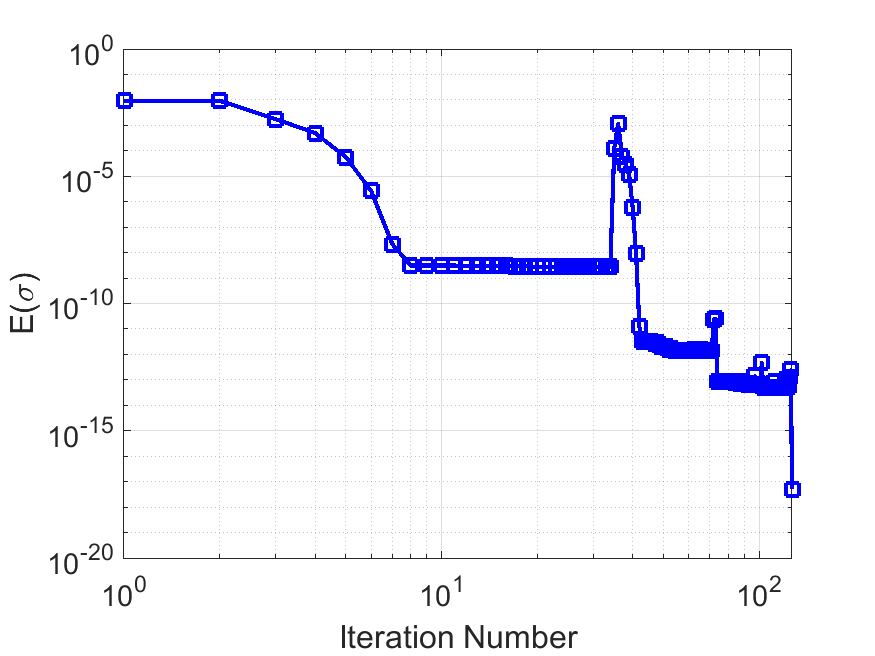}
		\caption{Error \textit{vs} iteration number.}		
	\end{subfigure}%
	~ 
	\begin{subfigure}[b]{0.4\textwidth}
		\includegraphics[width=\columnwidth]{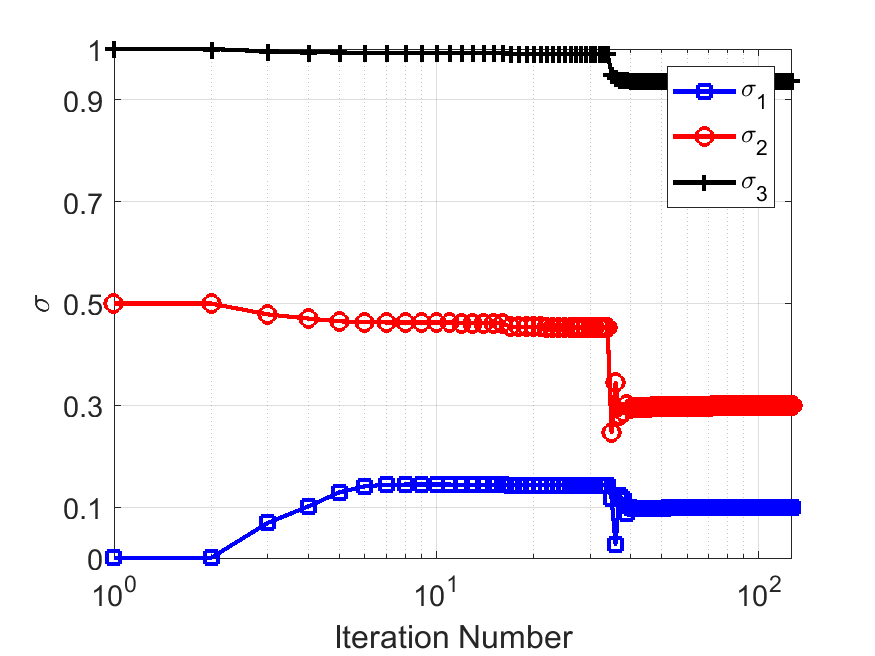}
		\caption{$\boldsymbol{\sigma}$ \textit{vs} iteration number.}		
	\end{subfigure}%
	\caption{Results for $M=3$, $S = 2$.  The obtained converged values are $\boldsymbol{\sigma} = \left[0.099997,\,\, 0.299995,\,\,0.939015\right]^T$, where $\sigma_1$, $\sigma_2$ captured the singularities, respectively with errors $E^\sigma_1 = 3.07\times 10^{-5}$ and $E^\sigma_2 = 1.65 \times 10^{-5}$. \label{fig:f12}}
\end{figure}

We test the case for $M=3$, using $S = 1,\,2,\,3$. In this analysis, we use $\tilde{N}=3$ initial data points. Figures \ref{fig:f11} and \ref{fig:f12} illustrate, respectively, the obtained results for $S=1,\,2$ using Algorithm \ref{alg:2}. We observe a quadratic convergence rate for $S = 1$, similar to Fig.\ref{fig:f5}, but with a few more iterations. Furthermore, for $M=2$, we observed a small change in the final value of $\sigma_3$, but nevertheless, capturing the two singularities $\sigma^*_1$, $\sigma^*_2$ is sufficient to minimize $E(\boldsymbol{\sigma})$ when $S = 2$.

\begin{figure}[t!]
	\centering
	\begin{subfigure}[b]{0.4\textwidth}
	\includegraphics[width=\columnwidth]{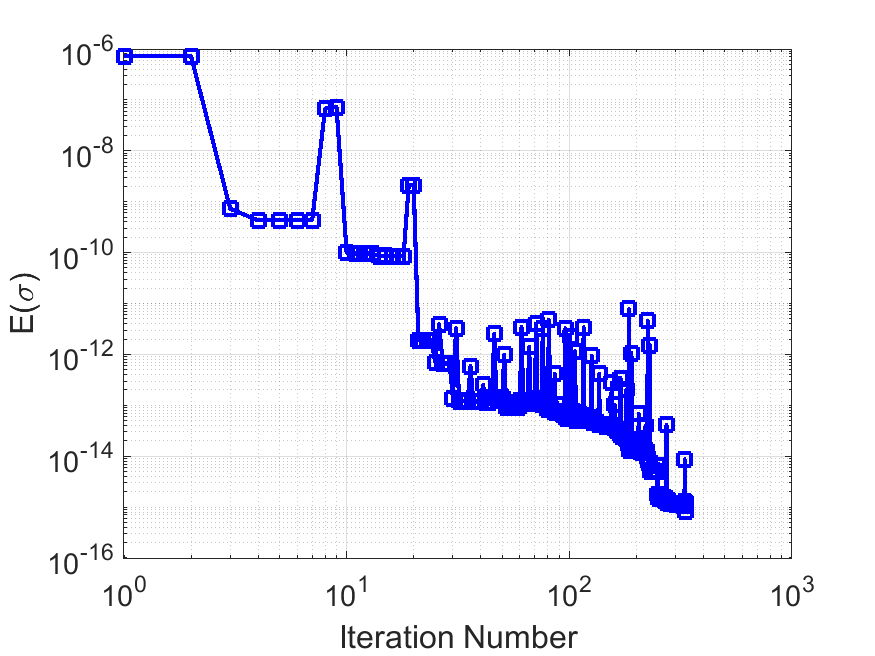}
		\caption{Error \textit{vs} iteration number.}		
	\end{subfigure}%
	~ 
	\begin{subfigure}[b]{0.4\textwidth}
		\includegraphics[width=\columnwidth]{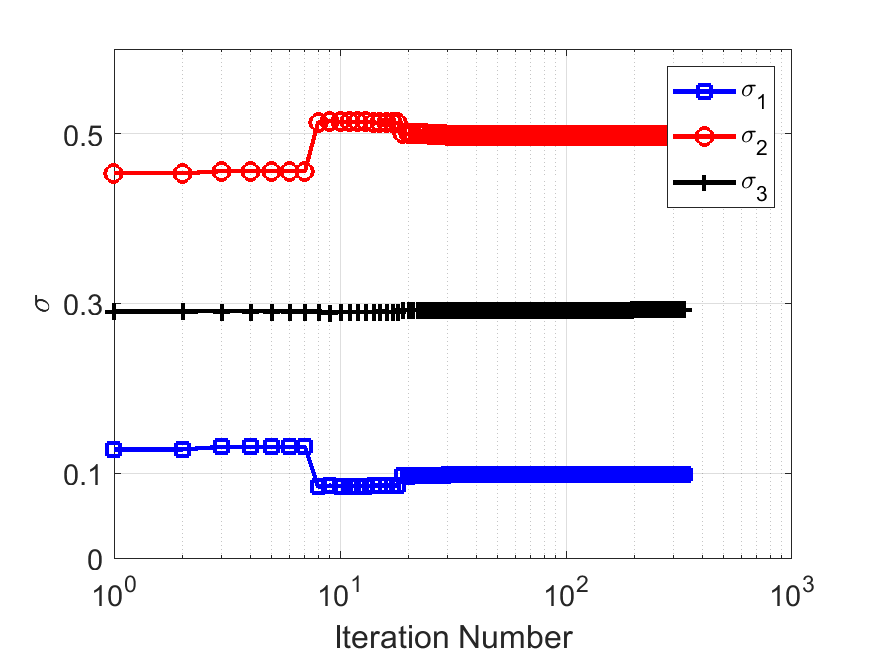}
		\caption{$\boldsymbol{\sigma}$ \textit{vs} iteration number.}		
	\end{subfigure}%
	\caption{Results for $M=3$, $S = 3$ using $\epsilon = 10^{-15}$. The converged values are $\boldsymbol{\sigma} = \left[0.099,\,\, 0.498,\,\, 0.293\right]^T$, which correspond to component-wise relative errors $E^\sigma = \left[0.014,\,\, 0.005,\,\, 0.024\right]^T$. \label{fig:f13}}
\end{figure}

The numerical convergence of Algorithm \ref{alg:2} towards the correct singularities $\boldsymbol{\sigma}^*$ becomes much more difficult when using $S = M = 3$. Therefore, we make use of the self-capturing approach summarized in Algorithm \ref{alg:1}, which incrementally solves the minimization problem using $M = 1, 2, 3$ and use the respective output singularities as initial guesses for the subsequent number of correction terms $M = M + 1$. We present the obtained results in Fig.\ref{fig:f13}. We observe that the scheme converges to the true singularities with errors for $\sigma$ up to $0.024$, but nevertheless, they yield an approximation error for $u^N$ of $10^{-15}$. 

\subsection{Random Singularities}
\label{Sec:Num_Random}

We perform two numerical tests involving random power-law singularities, where we define three strong singularities $\sigma^*_1$, $\sigma^*_2$, $\sigma^*_3$ randomly sampled from a uniform distribution $\mathcal{U}(0,\,1/2)$. We employ the self-singularity-capturing procedure in the context of single- and multi- term FDEs.

\subsubsection{FDEs with Random Singularities}
\label{Sec:random_1termFODE}

\begin{figure}[t!]
	\centering
	\includegraphics[width=0.5\textwidth]{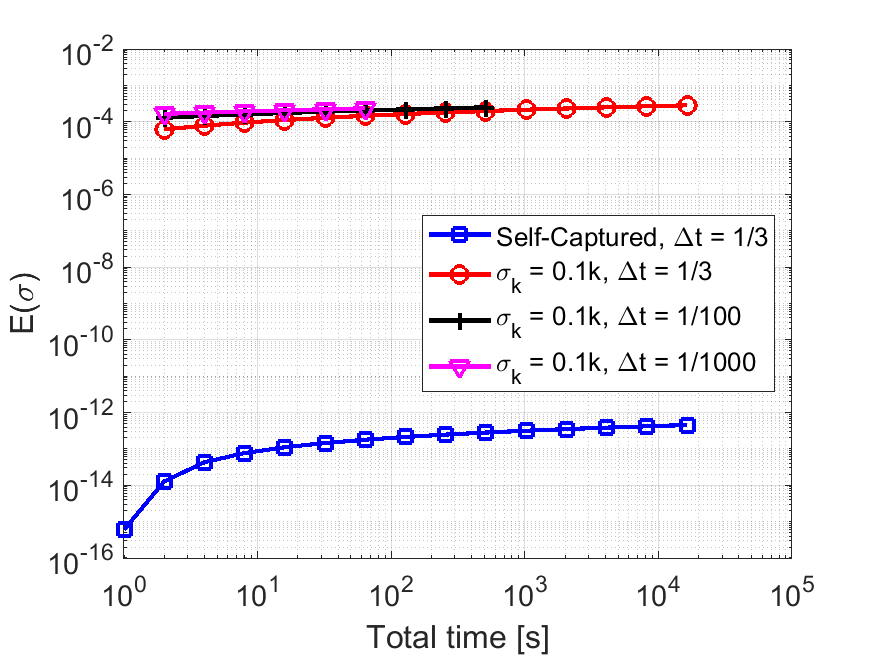}
	\caption{Comparison of error \textit{vs} time for long time-integration between the captured singularities and $\sigma_k = 0.1k$ for different values of time-step size $\Delta t$. \label{fig:f14}}
\end{figure}

We test the two-stage framework for efficient time-integration of (\ref{eq:problemdef}). We then compute $u^{ext}(t)$ and $f^{ext}(t)$, respectively using (\ref{eq:udelta}) and (\ref{eq:fdelta}), and use the self-capturing framework presented in Algorithm \ref{alg:1}. We then compare the approximation error $E(\boldsymbol{\sigma})$ with the choice $\sigma_k = 0.1 k$ defined in \cite{Zheng2017} when the singularities are unknown. Although our framework has no explicit information about the generated random singularities, we present them for verification purposes, which are given by:
\begin{equation}
\boldsymbol{\sigma}^* = \left[0.0172230402514543, \,\,\, 0.219372179828199, \,\,\, 0.190779228546504 \right]^T. 
\end{equation}
We set $\tilde{N}=3$ with $\Delta t = 1/3$, and therefore, $\tilde{\Omega} = [0, 1]$, $\epsilon = 10^{-15}$, and $\epsilon_1 = 10^{-13}$ to Algorithm \ref{alg:3}. Stage-I of the framework (Algorithm \ref{alg:1}) converges with $E(\boldsymbol{\sigma}^{(k)}) = 2.96 \times 10^{-16}$ with $M=2$ correction terms, with obtained values:
\begin{equation}\label{eq:captured}
 \boldsymbol{\sigma} = \left[0.0187990387914248, \,\,\,  0.206944449676742\right]^T. 
\end{equation}
After capturing the singularities within the desired precision $\epsilon$, we enter Stage-II and compute the numerical solution $u^N(t)$ using the captured $\boldsymbol{\sigma}$ for multiple, longer time-integration domains $\Omega$. We remark that when using an \textit{ad-hoc} choice of fixed $\sigma_k = 0.1 k$ with 4 correction terms, that is,
\begin{equation}\label{eq:power_fanhai}
\boldsymbol{\sigma} = \left[0.1, \,\,\, 0.2, \,\,\, 0.3, \,\,\, 0.4 \right]^T,
\end{equation}
we obtain an approximation error $E(\boldsymbol{\sigma}) = 5.25 \times 10^{-5}$ over $N = 4$ time-steps. We use the captured powers (\ref{eq:captured}) and the predefined ones (\ref{eq:power_fanhai}) for longer time-integration. The obtained results are presented in Figure \ref{fig:f14}, where we observe how using 3 time-steps to capture the singularities leads to precise long time-integration with a large time-step size $\Delta t = 1/3$ (blue curve). On the other hand, the errors for $\sigma_k = 0.1k$ are much larger, and also do not improve with smaller $\Delta t$ due to the presence of very strong singularities.

\subsubsection{Multi-Term FDEs with Random Singularities}
\label{Sec:multi-term}

In this section we show how the scheme in Stage-I can be applied for multi-term fractional differential equations. Let the following multi-term FDE:
\begin{equation}\label{eq:multitermFODE}
\sum^{N_\alpha}_{l=1} \prescript{RL}{0}{}\mathcal{D}^{\alpha_l}_t u(t) = f(t), \quad u(0) = u_0, \quad 0 < \alpha_l < 1.
\end{equation}
where $\alpha_l$ denotes the $N_\alpha$ multi-fractional orders with $l=1, \dots N_\alpha$. Using the approximation (\ref{eq:corrected_deriv}), and assuming the same powers $\boldsymbol{\sigma}$ for the correction weights for all fractional derivatives, we obtain:
\begin{equation}\label{eq:multitermdiscrete}
\sum^{N_\alpha}_{l=1} \prescript{RL}{0}{}\mathcal{D}^{\alpha_l}_t u^N(t) \big|_{t = t_n} + \sum^M_{j=1} \left( \sum^{N_\alpha}_{l=1} W^{(\alpha_l)}_{j,n}(\boldsymbol{\sigma}) \right) \left(u^N_j - u_0 \right) = f^ \delta_n,
\end{equation}
where we use the superscript $(\alpha_l)$ in $W^{(\alpha_l)}_{j,n}$ to denote the distinct sets of correction weights due to the fractional order $\alpha_l$ through (\ref{eq:final_weights}). The entire minimization scheme is identical, since it depends on the solution data for $u(t)$. In order to solve (\ref{eq:multitermdiscrete}) with the developed scheme in Stage-I, due to the linearity of (\ref{eq:multitermFODE}), we only need to replace the initial correction weights $W_{n,j}$, and $\alpha$-dependent coefficients $d_j,\,b_j$ in (\ref{eq:u13}) and (\ref{eq:ug3}) with the following summations:
\begin{equation}
\tilde{W}_{n,j} = \sum^{N_\alpha}_{l=1} W^{(\alpha_l)}_{n,j}, \quad \tilde{d}^{(p)}_j = \sum^{N_\alpha}_{l=1} d^{(p, \alpha_l)}_j, \quad
\tilde{b}^{(k)}_j = \sum^{N_\alpha}_{l=1} b^{(k, \alpha_l)}_j, \quad k = 1,\,2,\,3.
\end{equation}

We utilize the fabricated solution (\ref{eq:udelta}) in (\ref{eq:multitermFODE}), which yields the following force term:
\begin{equation}\label{eq:fdeltamulti}
f^{ext}(t) = \sum^{N_\alpha}_{\alpha_l} \sum^{S}_{j=1} \frac{\Gamma(1+\sigma^*_j)}{\Gamma(1+\sigma^*_j - \alpha_l)} t^{\sigma^*_j - \alpha_l}.
\end{equation}

We set $N_\alpha = 3$ fractional orders $\lbrace \alpha_1, \alpha_2, \alpha_3 \rbrace = \lbrace 0.3, 0.5, 0.7 \rbrace$, and identically to Section \ref{Sec:random_1termFODE}, we sample the following random singularities from $\mathcal{U}(0,1/2)$:
\begin{equation}
\boldsymbol{\sigma}^* = \left[0.13924910943352420, \,\,\, 0.2734407596024919, \,\,\, 0.4787534177171488 \right]^T. 
\end{equation}
Also, we set $\tilde{N}=3$ with $\Delta t = 1/3$ and $\tilde{\Omega} = [0, 1]$, $\epsilon = 5 \times 10^{-15}$, and $\epsilon_1 = 10^{-14}$, with $\gamma^0 = 10^ {-3}$ to Algorithm \ref{alg:1}, where the scheme converges with the following approximate singularities:
\begin{equation}
\boldsymbol{\sigma} = \left[0.1469249923105880, \,\,\, 0.4869203803691072, \,\,\, 0.3066386453671829 \right]^T. 
\end{equation}

\begin{figure}[t!]
	\centering
	\begin{subfigure}[b]{0.325\textwidth}
	\includegraphics[width=\columnwidth]{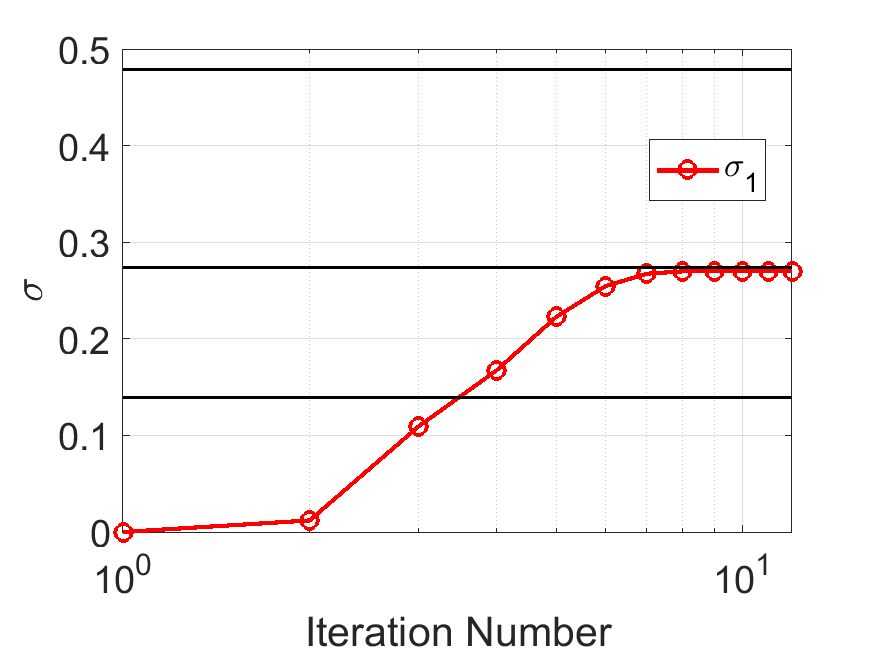}
		\caption{$M=1$}		
	\end{subfigure}
	\begin{subfigure}[b]{0.325\textwidth}
		\includegraphics[width=\columnwidth]{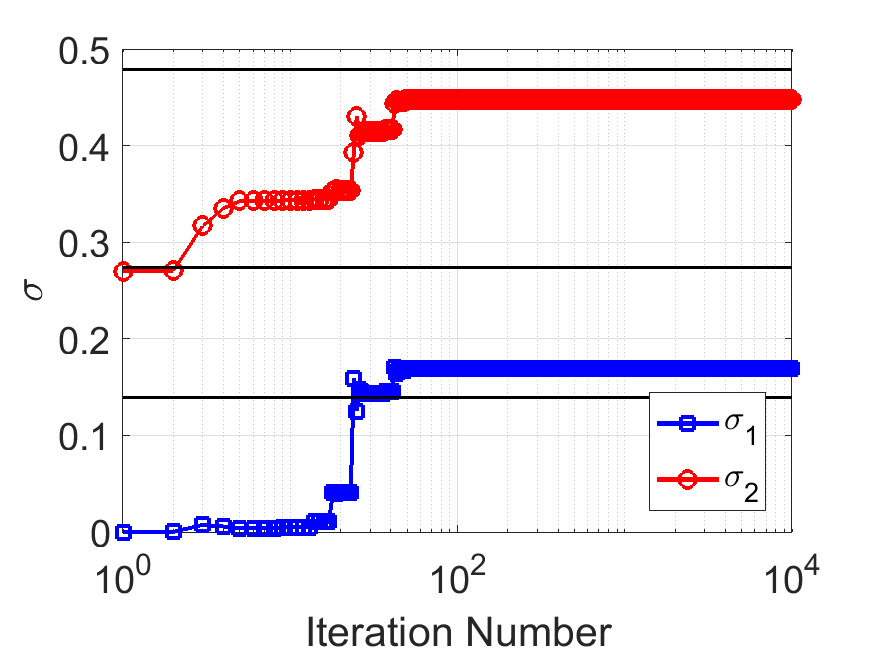}
		\caption{$M=2$}		
	\end{subfigure}
	\begin{subfigure}[b]{0.325\textwidth}
		\includegraphics[width=\columnwidth]{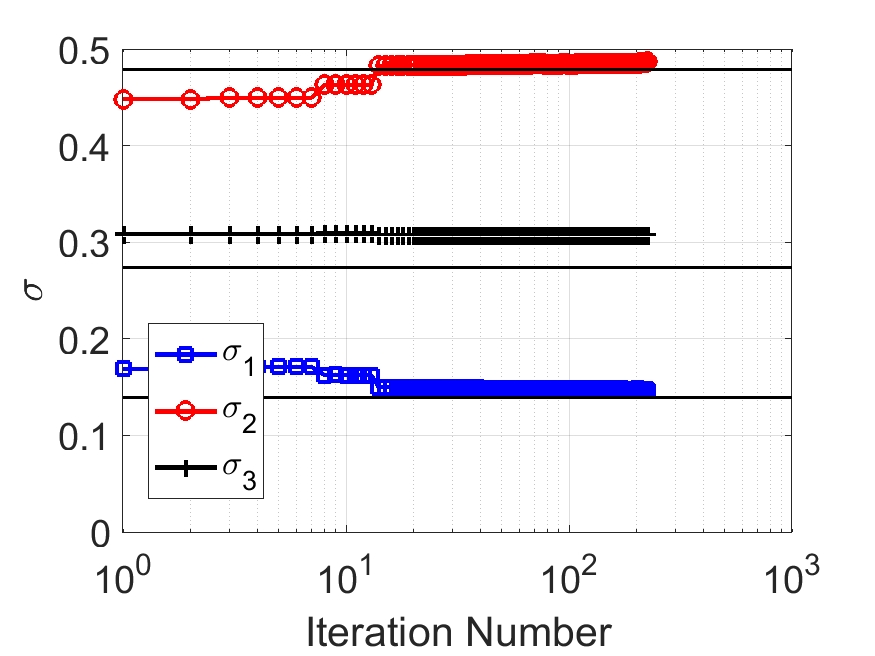}
			\caption{$M=3$.}		
		\end{subfigure}
	\caption{Singularities $\boldsymbol{\sigma}$ \textit{vs} iteration number for a multi-term FDE with $N_\alpha = 3$ in Stage-I for each $M$ in Algorithm \ref{alg:1}. \textit{(a)} $M=1$, with error $E(\boldsymbol{\sigma}) = 1.10 \times 10^{-5}$, \textit{(b)} $M=2$, with error $E(\boldsymbol{\sigma}) = 3.83 \times 10^{-13}$, and \textit{(c)} $M=3$ with convergence at $E(\boldsymbol{\sigma}) = 2.31 \times 10^{-15}$. The black horizontal lines correspond to the true random singularities in the solution. \label{fig:multiterm_random}}
\end{figure}

Figure \ref{fig:multiterm_random} presents the obtained results for each $M$ in the self-capturing procedure, where we observe that the scheme does not fully capture the true singularities, but provides sufficiently good approximations to obtain errors as low as $E(\boldsymbol{\sigma}) = 2.31 \times 10^{-15}$.

\subsection{Singular-Oscillatory Solutions}
\label{Sec:Oscillatory}

Finally, we test the two-stage framework considering a fabricated solution for (\ref{eq:problemdef}) to be the following multiplicative coupling between a power-law and oscillatory parts:
\begin{equation}\label{eq:oscillatory}
u^{ext}(t) = t^{\sigma^*} \cos(\omega t),
\end{equation}
where the corresponding right-hand-side is analytically obtained as:
\begin{equation}\label{eq:rhs_oscillatory}
f^{ext}(t) = c_1 t^{\sigma^*-\alpha} \left[{}_{p}\tilde{F}_q (a^{(1)}; b^{(1)}; -\omega^2 t^2/4) + c_3 \omega^2 t^2  {}_{p}\tilde{F}_q (a^{(2)}; b^{(2)}; -\omega^2 t^2/4) \right],
\end{equation}
with $c_1 = -2^{\alpha-\sigma^*-4}\pi \Gamma(1+\sigma^*)$, $c_2 = 8\left(\alpha - \sigma^* - 1\right)$ and $c_3 = (1+\sigma^*)(2+\sigma^*)$. Also, ${}_{p}\tilde{F}_q$ represents the regularized, generalized hypergeometric function given by \cite{Weisstein2003}:
\begin{equation}
{}_{p}\tilde{F}_q(a_1, \dots, a_p; b_1, \dots, b_q; z) = \frac{{}_{p}F_q(a_1, \dots, a_p; b_1, \dots, b_q; z)}{\Gamma(b_1) \dots \Gamma(b_q)},
\end{equation}
where, for the particular case (\ref{eq:rhs_oscillatory}), we have:
\begin{equation} 
	\begin{cases}
		a^{(1)} = \lbrace (1+\sigma^*)/2,\, (2+\sigma^*)/2 \rbrace, &
		b^{(1)} = \lbrace 1/2,\, (2-\alpha+\sigma^*)/2,\,(3-\alpha+\sigma^*)/2 \rbrace, \nonumber\\ \vspace{-2mm} \\
		a^{(2)} = \lbrace (3+\sigma^*)/2,\, (4+\sigma^*)/2 \rbrace, &
		b^{(2)} = \lbrace 3/2,\, (4-\alpha+\sigma^*)/2,\,(3-\alpha+\sigma^*)/2 \rbrace. \nonumber
	\end{cases}
\end{equation}

We fix the fractional order $\alpha = 0.5$, frequency $\omega = 10\pi$ and randomly sample $\sigma^*$ from a uniform distribution $\mathcal{U}(0,\,1/3)$ and obtain a value $\sigma^* = 0.2426481954401539$. For Stage-I, we consider $\tilde{N} = 3$ data points with $\tilde{\Omega} = [0, 0.01]$, with $\Delta t = 3/100\,[s]$. Furthermore, we set $\epsilon = 10^{-11}$, $\epsilon_1 = 10^{-14}$ and $\gamma^0 = 10^{-3}$. Using Algorithm \ref{alg:1} with a slight initial guess modification for $\boldsymbol{\sigma}^{(0)}$ when $M=2$ in line 11 to $\boldsymbol{\sigma}^{(0)} = [1,\,\sigma^{(k)}_1]^T$ (since the developed framework is mainly oriented to strong singularities), we capture the following two singularities:
\begin{equation}
\boldsymbol{\sigma} = \left[0.2427452349772425, \,\,\, 2.220682758797950 \right]^T, \nonumber
\end{equation}
where the first converged singularity is an approximation to the randomly sampled value for $\sigma^*$. To understand the second captured singularity $\sigma_2$, let the series expansion of (\ref{eq:oscillatory}) around $\omega t$:
\begin{equation}\label{eq:exact_oscillatory}
u^{ext}(t) = t^{\sigma^*} \cos(\omega t) = t^{\sigma^*} \left(1-\frac{\omega^2}{2} t^2 + \mathcal{O}(t^4) \right) = t^{\sigma^*} -\frac{\omega^2}{2} t^{2+\sigma^*} + \mathcal{O}(t^{4+\sigma^*}).
\end{equation}
Therefore, the captured $\sigma_2$ is an approximation of the second term $2 + \sigma^*$ of the above series expansion. Although capturing the third singularity $4+\sigma^*$ is possible, we would need to set higher values for the initial guesses in the self-capturing algorithm, and finding the corresponding minimum becomes more difficult, since such higher-order, weaker singularities have less influence in the accuracy of the solution. 

\begin{figure}[t!]
	\centering
	\begin{subfigure}[b]{0.4\textwidth}
	\includegraphics[width=\columnwidth]{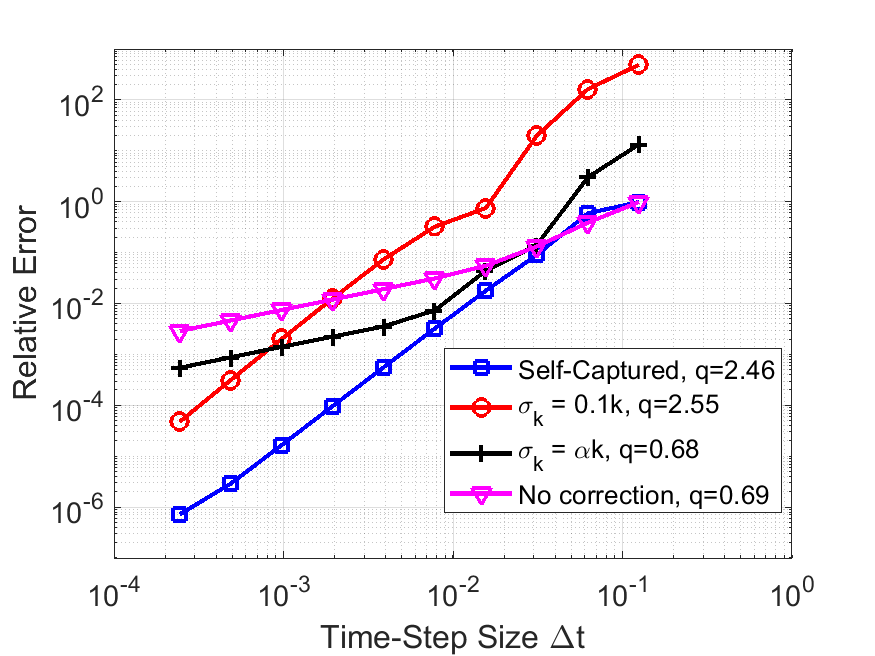}
		\caption{Relative error \textit{vs} time-step size.\label{fig:conv_oscillatory}}		
	\end{subfigure}%
	~ 
	\begin{subfigure}[b]{0.4\textwidth}
		\includegraphics[width=\columnwidth]{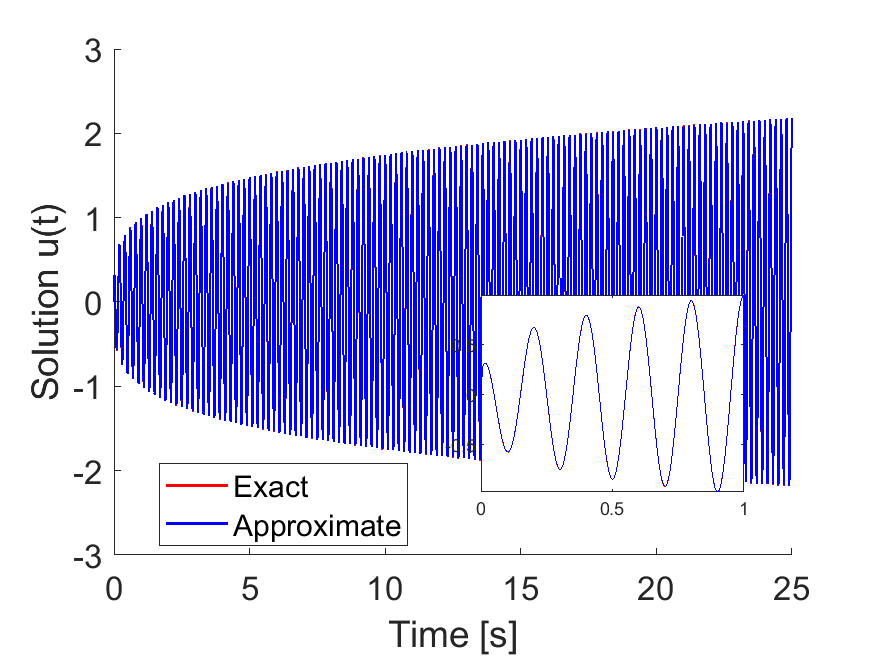}
		\caption{Solutions for $u(t)$ \textit{vs} time.\label{fig:oscillatory_longerint}}		
	\end{subfigure}%
	\caption{\textit{(a)} Convergence analysis comparing the accuracy of the captured singularities in Stage-I with \textit{ad-hoc} choices for time-integration with $T=1\,[s]$ in Stage-II. The term $q$ denotes the convergence rate. The self-captured case uses $M=2$ correction terms, while for $\sigma_k = 0.1k$ and $\sigma_k = \alpha k$ we use $M=4$. \textit{(b)} Exact and approximate solutions for $u(t)$ using the captured singularities for longer time-integration over $T=25\,[s]$. \label{fig:oscillatory}}
\end{figure}

After capturing the singularities in Stage-I, we apply them in Stage-II for a convergence test with $T = 1\,[s]$, with $\Omega = [0,T]$, and compare the approximate solution $u^N(t)$ with the exact one (\ref{eq:exact_oscillatory}) through the relative error defined as
\begin{equation}\label{eq:rel_error}
E_{L^2} = \frac{|| \mathbf{u}^N - \mathbf{u}^{ext} ||_{L^2(\Omega)}}{|| \mathbf{u}^{ext} ||_{L^2(\Omega)}}.
\end{equation}
between the two-stage framework other \textit{ad-hoc} choices for $\boldsymbol{\sigma}$, and also without corrections (see Figure \ref{fig:conv_oscillatory}). We observe that capturing $M=2$ singularities leads to the theoretical accuracy $\mathcal{O}(\Delta t^{3-\alpha})$ of the employed discretization, and also to errors about 2 orders of magnitude lower than the choice $\sigma_k = 0.1k$ with $M=4$ for the singularities in the observe range for $\Delta t$. We also compare the exact and approximate solutions over $T=25\, [s]$ with $N = 4096$ time-steps of size $\Delta t \approx 6.10 \times 10^{-3}$, and observe that we fully capture the singular and oscillatory behaviors of the solution.

\section{Conclusions}
\label{Sec:Conclusions}

We developed a two-stage framework for accurate and efficient solution of FDEs. It consists of: \textbf{I)} A self-singularity-capturing scheme, where the inputs were limited data for diminutive time, and the output was the captured power-law singular values of the solution. We developed an implicit finite-difference algorithm for the FDE solution and employed a gradient optimization scheme. \textbf{II)} The captured singularities were then inputs of an implicit finite-difference scheme for FDEs with accuracy $\mathcal{O}(\Delta t^{3-\alpha})$. Based on our numerical tests, we observed that:
\begin{itemize}
\item The scheme was able to fully capture $S$ singularities using $M \ge S$ correction terms. When $S > M$, the scheme still obtained approximations of the $M$ singularities, which was useful for error control.

\item When $M=1$, an estimation of the most critical singularity was captured by the scheme, regardless of $S$.

\item The developed scheme was used in a self-singularity-capturing framework, which determined the singularity values up to the desired tolerance, here set as $\epsilon = 10^{-15}$.

\item Once the singularities are captured in the first stage, we can successfully employ any known numerical scheme for more accurate and efficient solution of the corresponding FDE, when compared to the \textit{ad-hoc} singularity choices done in \cite{Zheng2017}.

\item The scheme was able to capture two singularities from a solution using a strong power-law singularity coupled with an oscillatory smooth part, which allowed us to obtain the theoretical accuracy $\mathcal{O}(\Delta t^{3-\alpha})$ of the employed numerical discretization for the fractional derivative.
\end{itemize}
The total computational complexity of the two-stage framework is $\mathcal{O}(\tilde{N}^2 N_{it} + N^2)$, since the number of correction terms $M$ is small. Among a variety of numerical schemes, Stage-II of the framework could incorporate fast convolution frameworks. This would dramatically increase the efficiency for long time-fractional integration of anomalous materials under complex loading conditions \cite{Suzuki2016, Suzuki2017Thesis} and also further improve the numerical accuracy for nonlinear problems using fractional multi-step schemes \cite{Zayernouri2016MS}. The developed scheme can also be used to capture the power-law attenuations and creep-ringing \cite{Jaishankar2013} in highly-oscillatory vibration of anomalous systems. Also, discovering the multi-singular behavior of the solution can give insight on knowledge- and kernel-based grid refinement near $t=0$ \cite{Stynes2017}. Furthermore, the presented scheme can leverage the development of self-constructing function spaces using M\"{u}ntz polynomials \cite{Esmaeili2011, Hou2017}, increasing the efficiency of spectral element methods \cite{kharazmi2017FSEM} applied to anomalous transport.

\appendix
\section{Discretization Coefficients for Fractional Derivatives}
\label{Ap:Discr_Coefs}

We present the $\alpha$- and $\Delta t$- dependent coefficients for the local (\ref{eq:discrlocal}) and history (\ref{eq:discrhist}) approximations for the RL fractional derivative, which are given by \cite{Zheng2017}:
\begin{equation} \label{eq:d2}
	\begin{cases}
		d^{(1)}_0 = -\frac{\alpha \Delta t^\alpha}{\Gamma(2+\alpha)}, \quad
		d^{(1)}_1 = \frac{\Delta t^\alpha}{\Gamma(2+\alpha)},		  														 & p=1, \\ \vspace{-2mm} \\
		d^{(2)}_0 = -\frac{\alpha \Delta t^\alpha}{2\Gamma(3+\alpha)}, \quad
	    d^{(2)}_1 = \frac{\alpha (3+\alpha) \Delta t^\alpha}{\Gamma(3+\alpha)}, \quad 
		d^{(2)}_2 = \frac{(4+\alpha) \Delta t^\alpha}{2\Gamma(3+\alpha)}, 
		& p=2.
	\end{cases}
\end{equation}
\begin{align}
b^{(1)}_j & = \frac{\Delta t^\alpha}{2 \Gamma(\alpha)}\left[a^{(\alpha+2)}_j - \left(2j-1\right)a^{(\alpha+1)}_j + j\left(j-1\right)a^{(\alpha)}_j \right], \label{eq:b1}\\
b^{(2)}_j & = -\frac{\Delta t^\alpha}{\Gamma(\alpha)}\left[a^{(\alpha+2)}_j - 2j a^{(\alpha+1)}_j + \left(j+1\right) \left(j-1\right)a^{(\alpha)}_j \right],\\
b^{(3)}_j & = \frac{\Delta t^\alpha}{2 \Gamma(\alpha)}\left[a^{(\alpha+2)}_j - \left(2j+1\right)a^{(\alpha+1)}_j + j\left(j+1\right)a^{(\alpha)}_j \right],\\
a^{(\alpha)}_j & = \frac{1}{\alpha} \left[\left(j+1\right)^\alpha - j^\alpha \right]. \label{eq:b3}
\end{align}

\section{Singularity-Capturing Formulation for a Single Correction Term and Time-Step}
\label{Ap:Single}

Let the time domain $\tilde{\Omega} = [0, \Delta t]$, and let $E:\mathbb{R} \to \mathbb{R}^+$ be the following quadratic error function over $\sigma$:
\begin{equation}
E(\sigma) = \vert\vert u^{data}(t) - u^N(t;\sigma) \vert\vert^2_{L^2(\tilde{\Omega)}},
\end{equation}
where we assume $u^{data}(t)$ to be known, and $u^N(t;\sigma)$ represents the numerical approximation of $u^{data}(t)$, with $\sigma \in \mathbb{R}$ to be determined, such that $E(\sigma)$ is minimized. 
Therefore, let the case for $M=1$ correction term and for the first time-step $t_1 = \Delta t$. We obtain the following cost function:
\begin{equation}\label{eq:cost_1step}
E(\sigma) = \left(u^{data}_1 - u^N_1(\sigma) \right)^2 = \left(\Delta t^{\sigma^*} - u^N_1(\sigma) \right)^2,
\end{equation}
where $u^N_1(\sigma)$ denotes the numerical solution for the FDE (\ref{eq:problemdef}) at $t = \Delta t$, obtained through (\ref{eq:discretized_initial}) and (\ref{eq:fdelta}) in the following fashion:
\begin{equation}\label{eq:init_system}
\prescript{RL}{0}{}\mathcal{D}^\alpha_t u^N(t)\big|_{t=\Delta t} +  W_{1,1} \left(u^N_1 - u_0\right) = \sum^{S}_{j=1} \frac{\Gamma(1+\sigma^*_j)}{\Gamma(1+\sigma^*_j - \alpha)}\Delta t^{\sigma^*_j - \alpha}.
\end{equation}
Recalling (\ref{eq:discretized_correction}), we obtain the initial correction weight $W_{1,1}$ by the following equation:
\begin{equation}
\prescript{RL}{0}{}\mathcal{D}^\alpha_t \left(t^\sigma\right)\big|_{t=\Delta t} + W_{1,1} \left(t^\sigma_1\right) = \frac{\Gamma(1+\sigma)}{\Gamma(1+\sigma-\alpha)} t^{\sigma-\alpha}_1,
\end{equation}
Recalling (\ref{eq:discrlocal}) with $p=1$:
\begin{equation}
d^{(1)}_1 \Delta t^\sigma + \Delta t^{\sigma} W_{1,1} = \frac{\Gamma(1+\sigma)}{\Gamma(1+\sigma-\alpha)} \Delta t^{\sigma-\alpha},
\end{equation}
therefore, we obtain the closed form for the correction weight:
\begin{equation}\label{eq:solution_W}
W_{1,1} = \frac{\Gamma(1+\sigma)}{\Gamma(1+\sigma-\alpha)}\Delta t^{-\alpha} - d^{(1)}_1,
\end{equation}
Substituting (\ref{eq:solution_W}) into (\ref{eq:init_system}), rewriting the RL derivative in its discretized form at $t=\Delta t$, and assuming homogeneous initial conditions $u_0 = 0$, we obtain:
\begin{equation}
d^{(1)}_1 u^N_1 + \left( \frac{\Gamma(1+\sigma)}{\Gamma(1+\sigma-\alpha)} \Delta t^{-\alpha} - d^{(1)}_1 \right) u^N_1 = f^{data}_1,
\end{equation}
hence,
\begin{equation} \label{eq:solu_1step}
u^N_1 = \Delta t^\alpha \frac{\Gamma(1+\sigma-\alpha)}{\Gamma(1+\sigma)} f^{data}_1.
\end{equation}
Substituting (\ref{eq:solu_1step}) into (\ref{eq:cost_1step}), we obtain, for $t=\Delta t$:
\begin{equation}\label{eq:cost1step}
E(\sigma) = \left(u^{data}_1 - \Delta t^\alpha \frac{\Gamma(1+\sigma-\alpha)}{\Gamma(1+\sigma)} f^{data}_1 \right)^2.
\end{equation}

We assume that the above equation has a root at $\sigma = \sigma^*$, that is, $E(\sigma^*) = 0$. Therefore, given an initial guess $\sigma^{(0)}$, we can apply a Newton scheme to iteratively obtain $\sigma^{k+1}$ in the following fashion:
\begin{equation}\label{eq:NRupdate}
\sigma^{k+1} = \sigma^k - \frac{E(\sigma)}{\frac{\partial E(\sigma)}{\partial \sigma}}\big|_{\sigma = \sigma^k},
\end{equation}
where $\frac{\partial E(\sigma)}{\partial \sigma}$ can be obtained analytically, and is given by:
\begin{equation}\label{eq:NRderivative}
\frac{\partial E(\sigma)}{\partial \sigma} = -2\psi^*\Delta t^\alpha \frac{\Gamma(1+\sigma-\alpha)}{\Gamma(1+\sigma)} \left(u^{data}_1 - \Delta t^\alpha \frac{\Gamma(1+\sigma-\alpha)}{\Gamma(1+\sigma)} \right),
\end{equation}
with $\psi^* = \psi_0(1+\sigma - \alpha) - \psi_0(1+\sigma)$, where $\psi_0(z)$ denotes the Polygamma function. We observe that the computational cost of the above procedure per iteration $k$ is minimal, since there is no history computation when evaluating the corresponding time-fractional derivatives of $u^N_1$.

\section*{Acknowledgments}

The authors would also like to thank Ehsan Kharazmi for productive discussions along the development of this work.

\bibliographystyle{siamplain}
\bibliography{SC_references}

\end{document}